\documentclass[english,a4paper,12pt]{article}

\usepackage[T1]{fontenc}
\usepackage[latin1]{inputenc}
\usepackage{babel,amsfonts,amsthm,epsfig,vmargin,amsmath,amssymb}
\usepackage{ae,aecompl}  

\setpapersize{A4}
\setmarginsrb{30mm}{14mm}{20mm}{15mm}{12pt}{9mm}{0pt}{5mm}

\newcommand{\R}{\mathbb{R}}

\newcommand{\K}{\mathbb{K}}
\newcommand{\F}{\mathbb{F}}

\newcommand{\Q}{\mathbb{Q}}

\theoremstyle{definition}
\newtheorem{thm}{Theorem} 
\newtheorem{lem}{Lemma}[section]
\newtheorem{define}{Definition}[section]
\newtheorem{coroll}{Corollary}[section]

\theoremstyle{remark}
\newtheorem{rem}{Remark}[section]

\begin{document}
\title{ Analysing singularities of a benchmark problem}
\author{{\sc Teijo Arponen},
  Institute of Mathematics, \\
  Helsinki University of Technology,
  PL 1100, 02015 TKK, Finland. \\  email: teijo.arponen@hut.fi. \\[8pt]
{\sc Samuli Piipponen},
  Department of Mathematics,
  University of Joensuu, \\
  PL 111, 80101 Joensuu, Finland. \\  email: samuli.piipponen@joensuu.fi. \\[8pt]
{\sc Jukka Tuomela},
  Department of Mathematics,
  University of Joensuu, \\
  PL 111, 80101 Joensuu, Finland. \\  email: jukka.tuomela@joensuu.fi.}
\pagestyle{headings}
\markboth{T. ARPONEN, S. PIIPPONEN, J. TUOMELA}{\rm Analysing singularities}
\date{\today}
\maketitle

\begin{abstract}
  The purpose of this paper is to analyze the singularities of a well
  known benchmark problem ``Andrews' squeezing mechanism''.  We show
  that for physically relevant parameter values this system admits
  singularities. The method is based on Gr\"obner bases computations
  and ideal decomposition. It is algorithmic and can thus be applied
  to study constraint singularities which arise in more general
  situations.
\end{abstract} 


{\bf Mathematics Subject Classification (AMS) 2000}: 70B15, 13P10, 70G25.

\section{Introduction} 

The ``Andrews' squeezing system'' was first described by Giles in
\cite{gi78:CTP} and further studied in \cite{ma81:CTS}.  It is a
planar multibody system whose topology consists of closed kinematic
loops (see Figure \ref{fig:angles}).  The Andrews' system was promoted
in \cite{sc90:MSH} as a benchmark problem to compare different
multibody solvers.  Nowadays it is a well-known benchmark problem
\cite{ha-wa91:SODE2,testset} for numerical integration of
differential-algebraic equations as well.  The equations are of the
Lagrangian form (or descriptor form, see also \cite{ar01:RCS})
\begin{equation}
  \label{eq:1}
\begin{cases}
  f(t,y,y',y'',\lambda) = 0 \\
  g(y) = 0
\end{cases}
\end{equation}
where the function $f$ describes the dynamical equations and $g$ gives
the (holonomic) constraints. Here $y\in\R^n$ are the (generalized) position
coordinates, $y'$ and $y''$ are the first and second derivatives,
respectively, and $\lambda$ is the Lagrange multiplier.

It is well known that singularities of any kind hinder solving
equations numerically \cite{ro-sc88:DMS,ha-wa91:SODE2,ba-av94:SFA,ei-ha95:RMC}.
Intuitively, a singularity is where the (generic) number of degrees of
freedom of the system changes. Mathematically these are the points where the rank of the Jacobian of $g$ drops. Hence in this paper we will not consider the actual dynamical equations and analyse only the constraints given by $g$.

Most differential equation solvers include a possibility to monitor
singularities, and usually when proximity of a singularity is
detected, the computation is best to be interrupted.  But this kind of
monitoring is local only, that is, it does not tell us a priori where
the singularities lie but only alert us when it is too late to fix
things, so to speak.  Also, the monitoring is often a non-negligible
part of computational cost.  Therefore, it would be highly useful to
know {\em a priori} where the singularities are, or to make sure that
there are no singularities, or perhaps even remove them (for the
latter approach, see \cite{ar01:RCS}).  Locating singularities has
been studied also in \cite{mc00:GDL}.  If we cannot avoid or remove
the singularities, at least knowing where they are encountered is
helpful (indeed, necessary) when planning the computation without
interruptions.  One can then tune the chosen integration algorithm
such that the disturbing effect of the singularities is diminished,
for example by compensating the singularity of the Kepler problem by a
local change of variables as in \cite{le-re05:SHD} within the
computation.  Further techniques on compensating singularities in
multibody systems are gathered and concisely compared in \cite{ba-av94:SFA} and \cite{ei-ha95:RMC}.

The paper is organized as follows: in the next Section we present the
situation in detail and formulate the constraint equations in polynomial form.  Section \ref{sec:algebra} gathers the necessary
algebraic tools.  Section \ref{sec:analysis} contains the actual
analysis where we show that the mechanism indeed has singularities for certain parameter values. In Section \ref{sec:examples} there are some numerical examples of singular
configurations, and in Section \ref{sec:discuss} we summarize and
discuss the results, and address possible future work.

\section{Andrews' squeezing mechanism}

\begin{figure}[htbp]
  \begin{center}
    \epsfig{file=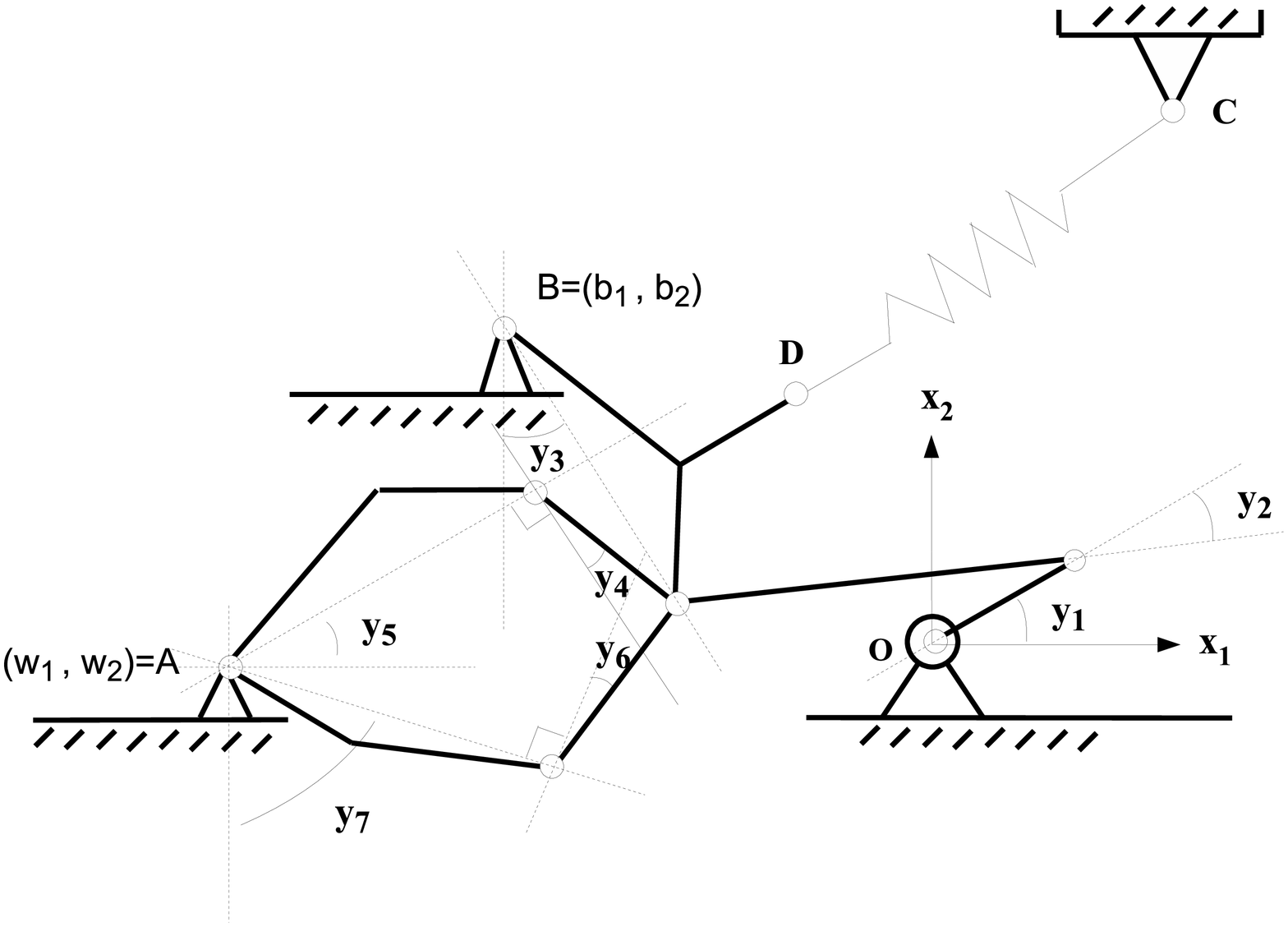,height=10cm,width=13cm}
  \end{center}
  \caption{The angles $y_i$ of the Andrews' system.}
  \label{fig:angles}
\end{figure}

\begin{figure}[htbp]
  \begin{center}
    \epsfig{file=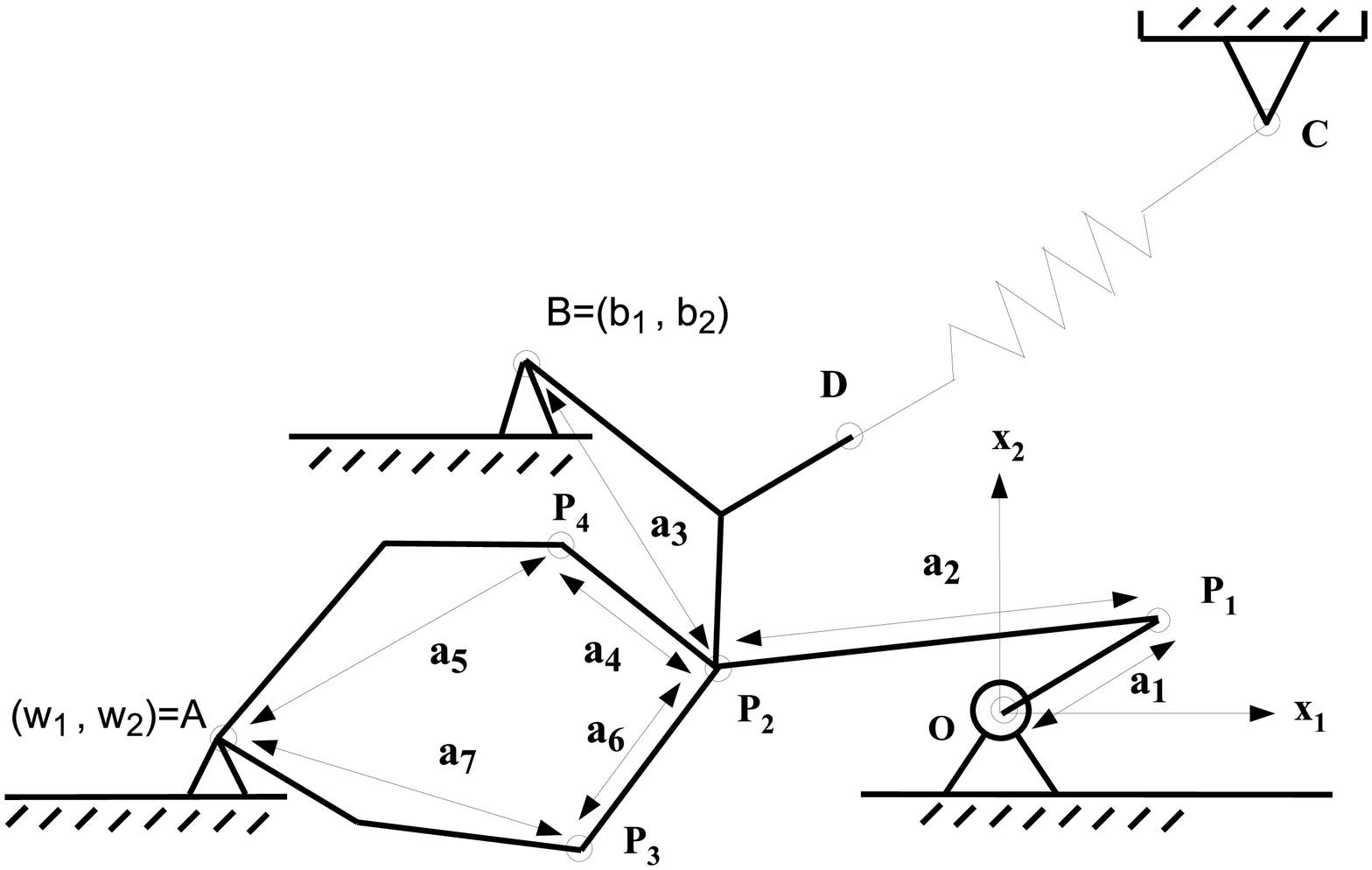,height=10cm,width=13cm}
  \end{center}
  \caption{The lengths $a_i$ and nodes of the Andrews' system.}
  \label{fig:lengths}
\end{figure}

The squeezing mechanism is given by the following equations.
\begin{equation}
  \label{eq:2}
  g(y) =
  \begin{cases}
    a_1\cos(y_1) - a_2\cos(y_1 + y_2) - a_3\sin(y_3) - b_1 \\
    a_1\sin(y_1) - a_2\sin(y_1 + y_2) + a_3\cos(y_3) - b_2 \\
    a_1\cos(y_1) - a_2\cos(y_1 + y_2) - a_4\sin(y_4 + y_5) 
    - a_5\cos(y_5) - w_1 \\
    a_1\sin(y_1) - a_2\sin(y_1 + y_2) + a_4\cos(y_4 + y_5) 
    - a_5\sin(y_5) - w_2 \\
    a_1\cos(y_1) - a_2\cos(y_1 + y_2) - a_6\cos(y_6 + y_7) 
    - a_7\sin(y_7) - w_1 \\
    a_1\sin(y_1) - a_2\sin(y_1 + y_2) - a_6\sin(y_6 + y_7) 
    + a_7\cos(y_7) - w_2
  \end{cases}
\end{equation}
Compared to the original articles mentioned above, we have chosen the
following notation for the parameters and angles:
\begin{align*}
  & a_1=rr\quad a_2=d  \quad a_3=ss\quad
  a_4=e  \quad a_5=zt\quad a_6=zf  \quad a_7=u  \\
  & b_1=xb\quad b_2=yb\quad  w_1=xa\quad w_2=ya \\
  & y_1 = \beta \quad y_2=\Theta \quad y_3=\gamma \quad
  y_4=\Phi \quad y_5=\delta \quad y_6=\Omega \quad y_7=\epsilon
\end{align*}
so the positions in Cartesian coordinates of the fixed nodes $A$ and
$B$ are given by $b=(b_1,b_2)$ and $w=(w_1,w_2)$, and the lengths of the rods by $a=(a_1,\dots,a_7)$, see Figures
\ref{fig:angles} and \ref{fig:lengths}.

Fixing the parameters $a$, $b$, and $w$, we have a map $g\,:\,\R^7\to\R^6$. Hence the set of possible  configurations, which is the zeroset $M_g=g^{-1}(0)$, is in general a curve (or possibly empty). Our task is to analyse the singularities of $M_g$, so let us state more precisely what is meant by a singularity.  
As mentioned before, in a singularity the number of degrees of freedom changes.  It is well known \cite{ro-sc88:DMS,ba-av94:SFA,mc00:GDL} that this corresponds to the situation where the rank of Jacobian drops.
\begin{define} Let $f\,:\,\R^n\to\R^k$ be any smooth map where $k<n$ and let $df$ be its Jacobian matrix.  Let $M=f^{-1}(0)\subset\R^n$ be the zeroset of $f$. A point $q\in M$ is a \emph{singular point} of $M$, if $df$ does not have maximal rank at $q$.
\end{define}
What in fact geometrically ``happens'' at a singular point may be quite complicated to determine. Typically the tangent space to $M$ does not change continuously in the neighbourhood of a singular point, or possibly $M$ intersects itself there. However, in all cases numerical problems occur, so it is important to try to find all singular points.

Note that the constraint equations \eqref{eq:2} (and hence the elements of its Jacobian matrix) are {\em not}
polynomials, yet our algebraic approach works only in a polynomial setting.
However, this problem is circumvented by reformulating $g(y)$ as
polynomials in the sines and cosines of $y_i$ by using the
trigonometric identities
\begin{align*}
  & \cos(x)^2+\sin(x)^2=1\\
  & \sin(x\pm y)=\sin(x)\cos(y)\pm\cos(x)\sin(y) \\
  & \cos(x\pm y)=\cos(x)\cos(y)\mp\sin(x)\sin(y)
\end{align*}
Setting $c_i=\cos(y_i), \quad s_i=\sin(y_i)$ we get the equations
\begin{equation}
  p(c,s)=\begin{cases}
   a_1  c_1-a_2\big(c_1c_2-s_1s_2\big)-a_3s_3-b_1=0\\
   a_1  s_1-a_2\big(s_1c_2+c_1s_2\big)+a_3c_3-b_2=0\\
   a_1  c_1-a_2\big(c_1c_2-s_1s_2\big)-a_4\big(s_4c_5+c_4s_5\big)
        -a_5c_5-w_1=0\\
   a_1  s_1-a_2\big(s_1c_2+c_1s_2\big)+a_4\big(c_4c_5-s_4s_5\big)
        -a_5s_5-w_2=0\\
   a_1  c_1-a_2\big(c_1c_2-s_1s_2\big)-a_6\big(c_6c_7-s_6s_7\big)
         -a_7s_7-w_1=0\\
   a_1  s_1-a_2\big(s_1c_2+c_1s_2\big)-a_6\big(s_6c_7+c_6s_7\big)
         +a_7c_7-w_2=0\\
    c_i^2+s_i^2-1=0,\qquad i=1,\dots,7.
\end{cases}
\label{psysteemi}
\end{equation}
We have 13 polynomial equations ($p_i=0$), 11 parameters
($a_1,\dots,a_7,\,b_1,b_2,w_1,w_2$) and 14 variables
($c_1,s_1,\dots,c_7,s_7$).  Note that each $p_i$ is of degree two in
$c_i,s_i$.  The equations $p_1=0,\dots,p_6=0$ correspond directly to
the 6 original equations $g(y)=0$ with the simple substitutions above (for
example $\cos(y_1+y_2)=c_1c_2-s_1s_2$) and the equations
$p_7=0,\dots,p_{13}=0$ are the extra identities due to ``forgetting''
the angle variables $y_i$.

Note that this reformulation of the constraints as algebraic equations is not just a trick which happens to work in this special case; indeed most constraints appearing in the simulation of multibody systems are of this type.

Now the above equations define $p$ as a map $p\,:\,\R^{14}\to\R^{13}$. Hence we expect that the zeroset $V=p^{-1}(0)\subset\R^{14}$ is a curve (or possibly empty). Singularities are then the points of this curve where the rank of $dp$ is not maximal. To find these points we need now to introduce some tools from commutative algebra.


\section{Background}
\label{sec:algebra}

In this section we present briefly the necessary definitions from
commutative algebra and algebraic geometry.  More details can be found
in \cite{cox-li-os92:IVA}, \cite{singularbook}, \cite{no76:FFR}, and
\cite{ei96:ca}.  These are roughly in the order of increasing difficulty,
\cite{cox-li-os92:IVA} being the most accessible, but unfortunately not containing  the necessary material on the Fitting ideals.

\subsection{Ideals and varieties}

Let $\K$ be an algebraic field and let $\K[x_1,\,\dots,\,x_n]$ be the ring of 
    polynomials in $x_1,\,\dots,\,x_n$, with coefficients in $\K$.
A subset $I \subset \K[x_1,\,\dots,\,x_n]$ is an
    {\em ideal} if it satisfies
    \begin{enumerate}
    \item[(i)] $0 \in I$.
    \item[(ii)] If $f,g \in I$, then $f+g \in I$.
    \item[(iii)] If $f \in I$ and $h \in
      \K[x_1,\,\dots,\,x_n]$, then $hf \in I$.
    \end{enumerate}
Ideals are often given by \emph{generators}. Let $f_1,\dots,f_s \in \K[x_1,\,\dots,\,x_n]$.  Then the set
    \begin{equation*}
      \langle f_1,\dots,f_s \rangle := \left\{ \sum_{i=1}^s h_i f_i 
        \mid  h_1,\dots,h_s \in \K[x_1,\,\dots,\,x_n] \right\}
    \end{equation*}
is an {\em ideal generated by} $f_1,\dots,f_s$. Any set of generators is called a \emph{basis}. 

Ideals are purely algebraic objects.  The geometrical counterpart of an
ideal is its locus, or variety. Let $I$ be an ideal in $ \K[x_1,\,\dots,\,x_n]$.
    Its corresponding {\em variety} is
    \[ 
    \mathsf{V}_{\F}(I) = \{ (a_1,\dots,a_n)\in\F^n \mid f(a_1,\dots,a_n)=0
    \quad \forall f\in I \}
    \]
where $\F$ is some field extension of $\K$. Note that it is often natural to choose $\F$ different from $\K$. If the field is clear from context we will sometimes write simply $\mathsf{V}(I)$.

Now different ideals may have the same variety. However, if one is interested mainly in the variety then it is useful to define
\[
  \sqrt{I}=\big\{f\in\K[x_1,\dots,x_n]\,|\, f^n\in I \textrm{\ \ for\ some\ \ }n\ge 1\big\}.
\]
If $I$ is an ideal, then $\sqrt{I}$ is the \emph{radical} of $I$; it
is the biggest ideal that has the same variety as $I$ and all ideals
having the same variety have the same radical.  Also, always $I\subset
\sqrt{I}$ and if $I=\sqrt{I}$ we say that $I$ is a {\em radical
  ideal}.
Some rudimentary properties among ideals and their varieties are in
the following
\begin{lem} Let $I$ and $J$ be ideals.  Then
  \begin{enumerate}
  \item $\mathsf{V}(I \cup J) = \mathsf{V}(I) \cap \mathsf{V}(J)$.
  \item $\mathsf{V}(I \cap J) = \mathsf{V}(I) \cup \mathsf{V}(J)$.
  \item $I \subset J$ if and only if $\mathsf{V}(I) \supset \mathsf{V}(J)$.
 \end{enumerate}
\end{lem}
Next we have to express the rank condition algebraically. To this end we need
\begin{define}
If $I= \langle f_1,\dots,f_s \rangle$, its {\em Fitting ideal}
    $F_I$ is the ideal generated by all maximal minors of the Jacobian
    matrix of $(f_1,\dots,f_s)$.\footnote{In general one can define Fitting ideals of minors of any given size. However, the above definition is sufficient for purposes of the present paper.} 
\end{define}
Now $\mathsf{V}(F_I)$ corresponds to the points where the rank is not maximal. However, the points are required also to be on $\mathsf{V}(I)$. Hence we conclude that the set of singular points, $S$, is given by
\[
           S=\mathsf{V}(I\cup F_I)
\]
In analysing varieties it is often helpful to decompose them to simpler parts. Similarly one may try to decompose a given ideal to simpler parts. This leads to following notions. 
\begin{define}
 A variety $V$ is {\em irreducible} if $V=V_1 \cup V_2$ implies
    $V=V_1$ or $V=V_2$.\\
   An ideal $I$ is  {\em prime}  if
    $f,g\in\K[x_1,\,\dots,\,x_n]$ and $fg\in I$ imply that either $f\in I$ or
    $g\in I$.
\end{define}
There is a very close connection between prime ideals and irreducible varieties. The precise nature of this depends on the chosen field. However, for our purposes the following is sufficient.
\begin{lem} If $I$ is prime, then $\mathsf{V}(I)$ is irreducible.\\
Any radical ideal can be written uniquely as a finite
    intersection of prime ideals,
    \begin{equation*}
      \sqrt{I} = I_1 \cap \cdots \cap I_r,
    \end{equation*}
    where $I_i \not\subset I_j$ for $i \neq j$.  
\end{lem}
This is known as the
    {\em prime decomposition of $\sqrt I$} and the $I_i$'s are called the
    minimal associated primes  of $I$. The above Lemma then immediately gives:
\begin{coroll}
\[
         \mathsf{V}(I) = \mathsf{V}(\sqrt{I}) = \mathsf{V}(I_1) \cup \cdots \cup \mathsf{V}(I_r),
\]
where all $\mathsf{V}(I_i)$ are irreducible.  
\end{coroll}   
Hence our strategy in analysing varieties is to compute the minimal associated primes of the relevant ideal, and then examine each irreducible component separately.

\subsection{Gr\"obner bases}

An essential thing is that all the operations above, especially
  finding the radical and the prime decomposition can be computed {\em
    algorithmically} using the given generators of $I$. To do this we need to compute special bases for ideals, called Gr\"obner bases. We will only briefly indicate the relevant ideas and refer to \cite{cox-li-os92:IVA} and \cite{singularbook} for more details.

First we need to introduce {\em monomial orderings}. All the algorithms handling the
ideals are based on some orderings among the terms of the generators
of the ideal.

Intuitively, an ordering $\succ$ is such that given a set of
monomials (e.g. terms of a given polynomial), $\succ$ puts them in
order of importance: given any two monomials
$x^{\alpha}:=x_1^{\alpha_1}\dots x_n^{\alpha_n}$ and $x^{\beta}$,
where $\alpha\neq\beta$ are different multi-indices, then either
$x^{\alpha}\succ x^{\beta}$ or $x^{\beta}\succ x^{\alpha}$. A common choice is to use \emph{degree reversed
lexicographic} ordering \cite{cox-li-os92:IVA}. In our
analysis we shall frequently need {\em product orders}, which are
formed as follows: if $\succ_A$ and $\succ_B$ are two orderings, we
shall divide the variables $x_i$ into two subsets, and use $\succ_A$
on the first subset and $\succ_B$ on the second.  This is indicated
with the following notation:
\[
    \K[(x_4,x_5,x_7),(x_1,x_2,x_3,x_6)].
\]
This is the same set as $\K[x_1,\dots,x_7]$ but now
the parenthesis indicate that we will use $\succ_A$ among the
variables $(x_4,x_5,x_7)$, and $\succ_B$ among the variables
$(x_1,x_2,x_3,x_6)$, and moreover all monomials where variables of the first group appear are always bigger than monomials where there are only variables of the second group. We will see later why this is useful.

Finally, the aforementioned Gr\"obner basis is a special kind of
generating set, with respect to some ordering.  Given any set of
generators and an ordering, the corresponding Gr\"obner basis exists and can be computed. The relevant algorithm is usually called the \emph{Buchberger algorithm}. The drawback of this algorithm is that it has a very high complexity in the worst case, and in practice the complexity depends quite much on the chosen ordering.\footnote{So far,
  no satisfactory theory of Gr\"obner basis complexity has been done.}

Anyway Gr\"obner bases have proved to be very useful in many different applications. Nowadays there exist many different implementations and improvements of the Buchberger algorithm. We chose to use the 
  well-known program {\sf Singular} \cite{GPS05}, \cite{singularbook} in all the computations in this paper.

\section{Analysing singularities}
\label{sec:analysis}

\subsection{Geometric description of the singularities}
Now getting back to our system \eqref{psysteemi} we see that we can
take the components of $p$ to be elements of $\Q(a,b,w)[c,s]$ where
$\Q(a,b,w)$ is the field of rational functions of $a$, $b$, and
$w$. Hence we have an ideal $J=\langle p_1,\dots,p_{13}\rangle\subset
\Q(a,b,w)[c,s]$ and the corresponding Fitting ideal $F_J$. On the
other hand we may view the ``parameters'' $a$, $b$, and $w$ also as
variables since they appear polynomially in the equations; hence we
could also consider $J\subset \Q[a,b,w,c,s]$. Taking this point of
view we can give an intuitive description of what kind of situations
we can expect.
\[
   \begin{cases}
      J\subset \Q[a,b,w,c,s]\\
      V_{\R}( J)\subset \R^{25}.
   \end{cases}
\]
In this way $V_{\R}( J)$ should be 12 dimensional (recall $J$ is
generated by 13 equations), i.e. a curve depending on 11
parameters. On the other hand if we fix parameters $a$, $b$, and $w$ we
get a curve in $\R^{14}$ which will be denoted by $V_{a,b,w}$. In the
same way we can view $ V_{\R}( J\cup F_J)$ as a variety in $\R^{25}$,
and fixing the parameters we get the singular points $V_{a,b,w}^S$.
Obviously $V_{a,b,w}^S\subset V_{a,b,w}\subset\R^{14}$.

Then what kind of variety should  $V_{\R}( J\cup F_J)$ be? Since the Jacobian of $p$ is of size $13\times 14$, \emph{generically} we expect to get 2 independent conditions in order the rank to drop.  That is, augmenting $J$ with $F_J$ should bring in 2 more equations. Hence we expect that $V_{\R}( J\cup F_J)$ is 10 dimensional; in other words we expect that if 11 parameters are chosen independently  then $V_{a,b,w}^S$ should be empty. On the other hand if a single condition among parameters is satisfied, then $V_{a,b,w}^S$ should consist of isolated points.

Further, if there are 2 conditions among parameters (i.e. 9 parameters freely chosen), then it would be possible that $V_{a,b,w}^S$ were one dimensional. But then our original constraint equations would be redundant, i.e. there would be more than one degree of freedom.

Below we will in fact observe that if a certain condition on parameters is satisfied, $V_{a,b,w}^S$ is indeed a finite set of points.

\subsection{Singular variety}

To study $V_{\R}( J\cup F_J)$ we could in principle 
use Gr\"obner basis theory in a straightforward manner. Let $G$ be the Gr\"obner basis of $J\cup F_J$ using the product order $\Q [(c,s),(a,b,w)]$. Let us denote by $g_1,\dots,g_r$ the elements of $G$ which do not depend on $c$ and $s$.
\begin{define} Let $S_J=\langle g_1,\dots,g_r\rangle$; then we say that $\mathsf{V}_{\R}(S_J)\subset\R^{11}$ is the singular variety associated to $J$.
\end{define}
It follows from the Gr\"obner basis theory that $V_{a,b,w}$ can have singularities \emph{only if} $(a,b,w)\in \mathsf{V}_{\R}(S_J)$. Hence theoretically, we could now find the singularities of the Andrews' system
in a straightforward manner by calculating the Gr\"obner basis of
$J\cup F_{J}$.  But this is an enormous task, due to $F_J$
being generated by high degree polynomials, not to mention including
the 11 parameters $a,b,w$. We could not get the solution in a finite time using our work station with 64GB memory.

Instead,
something else needs to be done.  Luckily there is another approach:
noting that $p_1,p_3,p_5$ have common terms, as well as $p_2,p_4,p_6$,
gives us motivation to study two subsystems.  One spanned by $p_5-p_3$
and $p_6-p_4$, the other one spanned by $p_5-p_1$ and $p_6-p_2$ (along
with the relevant trigonometric identities from $p_7,\dots,p_{13}$).
These subsystems are handleable and give useful information for the
whole system as well. Proceeding in this way we could at least determine that the singular variety is not empty and we could compute some subvarieties of it.

\subsection{Subsystem 4567}
\label{sec:subsystem4567}

Intuitively, the nodes and bars 4, 5, 6, 7 formulate a subsystem, see Figures \ref{fig:angles} and \ref{fig:lengths}. We suspect that when the lengths $a_4,\dots,a_7$ are such that the ``4567''
system is able to become one-dimensional, hence in some sense
degenerated, there should be a singularity in the whole system (see also the net example
in \cite{ar01:RCS}).  We will
shortly see that this is indeed the case.

Define
\begin{align*}
  q_1 &:= p_5-p_3 
  = a_4\big(s_4c_5+c_4s_5\big)+a_5c_5-a_6\big(c_6c_7-s_6s_7\big)-a_7s_7 \\
  q_2 &:= p_4-p_6
  = a_4\big(c_4c_5-s_4s_5\big)-a_5s_5+a_6\big(s_6c_7+c_6s_7\big)-a_7c_7 \\
  q_{i} &:= p_{i+7}= c_{i+1}^2+s_{i+1}^2-1, \quad i=3,\dots,6.
\end{align*}
Note that $q_1,q_2$ contain only angles $c_i$, $s_i$ and parameters $a_i$ for $i=4,\dots,7$.  That is why we do not need the other $p_i$'s.  Let
$J_{4567}$ be the ideal spanned by $q_1,\dots,q_6$. Hence we have
\begin{equation}
  \label{eq:S1}
  J_{4567} \subset \Q[(c_4,s_4,c_5,s_5,c_6,s_6,c_7,s_7),(a_4,a_5,a_6,a_7)]
\end{equation}
where we have indicated the relevant product order. 
 The
Gr\"obner basis $G$ for $J_{4567}\cup F_{J_{4567}}$ with respect to this
ordering contains 191 elements (denoted by $g_1,\dots,g_{191}$), out
of which 3 are especially enlightening:
\begin{align*}
  g_5 &= c_6a_6a_7, \\
  g_{16} &= c_4a_4a_5, \qquad \text{ and} \\
  g_1 &= \prod_{i=1}^{8}t_i, \qquad \text{ where} \\
  &t_1=a_4-a_5-a_6-a_7\\
  &t_2=a_4-a_5+a_6+a_7\\
  &t_3=a_4+a_5+a_6+a_7\\
  &t_4=a_4+a_5-a_6-a_7\\
  &t_5=a_4-a_5+a_6-a_7\\
  &t_6=a_4-a_5-a_6+a_7\\
  &t_7=a_4+a_5-a_6+a_7\\
  &t_8=a_4+a_5+a_6-a_7.
\end{align*}
Since $g_1$ is the only generator which does not contain any variables $c_i$ and $s_i$ we conclude that
\begin{thm}
\label{thm:t_i}
 The singular variety of $J_{4567}$ is
\[
      S_{J_{4567}}=\mathsf{V}(\langle g_1\rangle).
\]
\end{thm}     
Note that the factorization of $g_1$ gives us the prime decomposition of $\langle g_1\rangle$ and hence decomposition of $\mathsf{V}(\langle g_1\rangle)$ into 8 linear irreducible varieties. 

Our next task is to show that at least some points of the singular variety extend to actual (physically relevant) singularities of the whole system. Recall that each generator $g_i$ corresponds to an equation $g_i=0$.
Since $a_i>0$ in physically relevant cases, generators $g_5$ and $g_{16}$ imply that all the singularities of $J_{4567}$
have necessarily $c_6=c_4=0$ (conditions for the angles $4$ and $6$). In other words, in ideal-theoretic language, we can as well study the
ideal
\[
T:=\langle J_{4567},\, F_{J_{4567}},\, c_4,\, c_6 \rangle.
\]
Now the prime decomposition of $\sqrt{T}$ has 16 components:
\begin{equation}
  \sqrt{T} = T_1 \cap \ldots \cap T_{16}.
\end{equation}
Inspecting the generators of each of $T_j$, it is noticed that every
$T_j$ contains the $t_i$'s or $a_i$'s.  Recall that a generator $a_i$
in an ideal corresponds in the variety to a condition $a_i=0$ which is
non-physical.  Moreover, $t_3$ is now a non-physical condition
contradicting $a_i>0\,\forall i$.  Hence we discard (as in
\cite{ar01:RCS}) those ideals which have a non-physical generator that
would imply $a_i\le 0$ for some $i$, and we are left with 7 ideals,
whose generators are:
\begin{align*}
  T_1 &=\langle c_7^2+s_7^2-1,\, t_1,\, s_6+1,\, s_5-c_7,\, c_5+s_7,\, s_4+1,\, c_4,\, c_6 \rangle \\
  T_2 &=\langle c_7^2+s_7^2-1,\, t_2,\, s_6+1,\, s_5+c_7,\, c_5-s_7,\, s_4+1,\, c_4,\, c_6 \rangle \\
  T_3 &=\langle c_7^2+s_7^2-1,\, t_4,\, s_6+1,\, s_5+c_7,\, c_5-s_7,\, s_4-1,\, c_4,\, c_6 \rangle \\
  T_4 &=\langle c_7^2+s_7^2-1,\, t_5,\, s_6-1,\, s_5-c_7,\, c_5+s_7,\, s_4+1,\, c_4,\, c_6 \rangle \\
  T_5 &=\langle c_7^2+s_7^2-1,\, t_6,\, s_6-1,\, s_5+c_7,\, c_5-s_7,\, s_4+1,\, c_4,\, c_6 \rangle \\
  T_6 &=\langle c_7^2+s_7^2-1,\, t_7,\, s_6-1,\, s_5-c_7,\, c_5+s_7,\, s_4-1,\, c_4,\, c_6 \rangle \\
  T_7 &=\langle c_7^2+s_7^2-1,\, t_8,\, s_6-1,\, s_5+c_7,\, c_5-s_7,\, s_4-1,\, c_4,\, c_6 \rangle.
\end{align*}
Especially, we see that $s_6=\pm 1$, $s_5=\pm c_7$, $c_5=\pm s_7$, and
$s_4=\pm 1$.  Now we are ready to continue with the original system
$J \cup F_J$. 
\begin{rem}\label{rem:T_i_physical} Mathematically speaking the analyses of all cases $T_i$ are completely similar. However, on physical grounds the cases $T_1$, $T_2$, $T_6$ and $T_7$ are not so interesting. Indeed, in these cases the length of one of the rods corresponding to $a_4$, $a_5$, $a_6$ and $a_7$ is equal to the sum of the lengths of three others. Hence all four rods could be modelled as a single rod which would make the whole model significantly simpler. In the remaining cases no such reduction can be done, and we chose to examine the ideal $T_5$ in detail. See also remark \ref{rem:all_Ti_same_Q}.
\end{rem}

The case $T_5$ gives us conditions $s_4=-1$, $s_6=1$, $s_5=-c_7$, $c_5=s_7$, and
$a_7=a_5+a_6-a_4$ which we substitute into the original system. Next we will show that the resulting system has real solutions. These will be the required singular points.

The above substitutions simplify the generators of $J\cup F_J$ so that we get the following ideal:
\begin{equation}
\label{kideaali}
\begin{aligned}
  K=&\ \langle K_1\cup K_2\rangle,\\
K_1\quad:&\quad\begin{cases}
  k_1=  a_2(-c_1c_2+s_1s_2)+c_1a_1-s_3a_3-b_1  \\
  k_2= a_2(-s_1c_2-c_1s_2)+s_1a_1+c_3a_3-b_2 \\
  k_3=  c_1^2+s_1^2-1  \\
  k_4= c_2^2+s_2^2-1,  
  \end{cases}\\
K_2\quad:&\quad\begin{cases}
  k_5=  s_7(a_4-a_5)+s_3a_3+b_1-w_1  \\
  k_6=  c_7(a_5-a_4)-c_3a_3+b_2-w_2  \\
  k_7=  c_3^2+s_3^2-1  \\
  k_8=  c_7^2+s_7^2-1.
  \end{cases} 
\end{aligned}
\end{equation}
In $K_2$ we have 4 equations for 4 unknowns $c_3$, $s_3$, $c_7$, and $s_7$; hence it appears reasonable that we can get a finite number of solutions. Then we can substitute the computed values to $K_1$ which then becomes also a system of 4 equations for 4 unknowns $c_1$, $s_1$, $c_2$, and $s_2$. By the same reasoning we again expect that it is possible to get some solutions for appropriate parameter values.

We could numerically solve the
variables from these equations (and, indeed, we will, in the numerical
examples), but to analyze the situation in more detail we need to study
these further.

Then starting with the system $K_2$ we solve the angles 3 and 7 by the
following trick.  First we inspect the ideal generated by $K_2$ in the
ring
\[
\Q(b_1,b_2,w_1,w_2,a_3,a_4,a_5)[c_3,s_3,c_7,s_7].
\]
Calculating the Gr\"obner basis $\tilde G$ of $\langle K_2\rangle$
with respect to the lexicographic ordering we get 4 generators:
\begin{equation}
\label{eq:gtilde}
\begin{aligned} 
  \tilde g_1&= f_1s_7^2+f_2s_7-f_3f_4 \\  
  \tilde g_2 &= 2(b_2-w_2)(a_4-a_5)c_7-2(b_1-w_1)(a_4-a_5)s_7+f_5=0\\
  \tilde g_3 &= a_3s_3+(a_4-a_5)s_7+b_1-w_1=0  \\
  \tilde g_4 &= a_3c_3+(a_4-a_5)c_7+w_2-b_2=0. 
\end{aligned}
\end{equation}
where the auxiliary expressions $f_i$ are lengthy combinations of the
parameters $a_i,b_i$ (see the appendix).\footnote{The algorithms
  actually give by default only sums of monomials instead of products like
  $2(b_2-w_2)(a_4-a_5)$ but we have simplified these by hand.  Also
  \textsf{Singular} \cite{GPS05} could be used to automatically
  factorize into products but would involve some more elaborate
  programming.}

Now $\tilde g_1$ contains only $s_7$ and parameters.  Note that $f_1=0$
if and only if $a_4=a_5$.   Assuming $a_4\neq a_5$
the equation $\tilde g_1=0$ is a polynomial in $s_7$ of degree 2, hence
in order to have real solutions we need to impose the condition
\begin{align}
\label{S1_T5_D_condition}
   f_2^2 +4 f_1 f_3 f_4 \geq 0.
\end{align}
This condition can easily be checked when the parameters $a,b,w$
have been given numerical values.  Once $s_7$ is known, $c_7,s_3,c_3$
can be solved from the linear equations of $\tilde G$, provided
$a_4\neq a_5$ and $w_2\neq b_2$.

The cases $w_2=b_2$ and/or $a_4=a_5$ can be summarized as follows:
\begin{itemize}
\item[(i)] If
$w_2=b_2$ but $a_4\neq a_5$, we still get equations similar to $\tilde
G$, but now $s_3$ has a quadratic equation instead of $s_7$. 
\item[(ii)] If $a_4=a_5$, the system typically does not have
  solutions.  At least, a further condition among parameters, namely
  $|b-w|=a_3$, arises.  We shall not elaborate this nongeneric
  behaviour further.  In Section \ref{sec:a4_is_a5} we consider an
  example of this situation.
\end{itemize}
\begin{rem}
  In general, when the inequality in \eqref{S1_T5_D_condition} is
  strict, $s_7$ has 2 possible values.  Therefore, the tuples
  $(s_3,c_3,s_7,c_7)$ have in general 2 possible values because the
  other ones in the tuple are determined uniquely from $s_7$.
\end{rem}
The only thing left to be done, in this $J_{4567}$ subsystem case, is to solve $c_1,s_1,c_2,s_2$.  This is done with the ideal $\langle K_1\rangle$ given in \eqref{kideaali}.
\begin{rem}
  \label{rem:all_Ti_same_Q}
  Had we used any other $T_i$ instead of $T_5$ above, we would have
  ended up with this same ideal $\langle K_1\rangle$.
\end{rem}
We calculate the Gr\"obner basis
$\hat G$ of $\langle K_1\rangle$ , this time in the ring
\[
\mathbb{Q}(a_1,a_2,a_3,b_1,b_2,c_3,s_3)[c_1,s_1,c_2,s_2].
\]
Note especially that $s_3,c_3$ are here treated as parameters, due to
being now known expressions in the parameters $a$, $b$, $w$.  We again
use lexicographic ordering and get 4 generators $\hat g_1,\dots,\hat g_4$. 
Analogously to $s_7$ above, now for $s_2$ we get the second degree
polynomial equation
\begin{equation}\label{eq:Q_s2}
  \hat g_{1}=(-4a_1^2a_2^2)s_2^2-n_1n_2=0
\end{equation} 
where
\begin{align*}
  n_1&=a_1^2+2a_1a_2+a_2^2 -a_3^2 -2a_3b_1s_3+2a_3b_2c_3-b_1^2-b_2^2\\
  n_2&=a_1^2-2a_1a_2+a_2^2 -a_3^2 -2a_3b_1s_3+2a_3b_2c_3-b_1^2-b_2^2
\end{align*}
and linear equations for $c_2,s_1,c_1$:
\begin{align*}
  &\hat g_2 = d_1c_2+d_2+d_3 \\
  &\hat g_3 = l_1s_1+l_2+l_3 \\
  &\hat g_4 = (a_1^2-a_2^2)c_1+ l_4 
\end{align*}
where the auxiliary expressions $d_i,\,l_i$ are certain known (but
lengthy) functions of $a,b$, apart from $l_4$ which depends on
$s_1,s_2,c_2$ as well.  (See the appendix.)  In order to 
have real solutions for $s_2$, \eqref{eq:Q_s2} implies the condition
\begin{align}
  \label{S1_T5_E_condition}
  E := n_1 n_2 \leq 0.
\end{align}
These $\hat g_i$ determine $s_2,c_2,s_1,c_1$ provided $d_1\neq 0$,
$l_1\neq 0$, $a_1\neq a_2$.  To analyse the cases $d_1=0$, $a_1=a_2$,
and/or $l_1=0$, it is helpful
to define
\[
d_0:=a_3^2 +2a_3b_1s_3 -2a_3b_2c_3 +b_1^2+b_2^2.
\]
It turns out that $l_1=0 \Leftrightarrow d_1=0 \Leftrightarrow
d_0=0$.  After rearranging the terms (see the appendix) it can be seen that
 the condition \eqref{S1_T5_E_condition}
is equivalent to
\[
 (a_1-a_2)^2 \le d_0 \le (a_1+a_2)^2.
\]
Therefore, if $a_1\neq a_2$ then $d_0\neq 0$ and the equations above
can be solved.  The case $a_1=a_2$, $d_0\neq 0$ does not essentially
change the situation: we still have a quadratic equation for $s_2$,
and linear ones for the others, with a different coefficient for
$c_1$.

The remaining case $a_1=a_2$, $d_0=0$ corresponds to the situation
where the centre node coincides with the origin.  This gives another
singularity (the angle $y_1$ remains arbitrary) but is a rather
special case and will not be pursued further here.

\begin{thm}
  \label{thm:1}
  Let us suppose that the parameters $a$, $b$, $w$ satisfy the following
  conditions: $a_4\neq a_5$ and
  \begin{align}
    & n_1 (4a_1a_2-n_1) \geq 0 \tag{\ref{S1_T5_E_condition}} \\
    & f_2^2 + 16(a_4-a_5)^2|b-w|^2 f_3f_4\geq 0  \tag{\ref{S1_T5_D_condition}}
  \end{align}
  Then $V_{a,b,w}$ contains at least 2 singular points. If the
  inequalities are strict we get in general at least 4 singular
  points.
\end{thm}
It may appear that we also have at most 4 singular points. However, it is a priori possible that the other systems $T_i$ yield more singular points with the same parameter values.
\begin{proof}
  The first part of the theorem merely collects what we have shown
  above, with the simplifications $n_2=n_1-4a_1a_2$ and
  $f_1=4(a_4-a_5)^2|b-w|^2$.  The conditions are due to 
  univariate second degree polynomial equations, which have real
  solutions if and only if \eqref{S1_T5_D_condition} and
  \eqref{S1_T5_E_condition} (for $s_7$ and $s_2$, respectively) are
  fulfilled.  The other variables are determined from linear
  equations: $s_4,c_4,\dots,s_6,c_6$ from $T_5$; $s_3,c_3,c_7$ from
  $K_1$; $s_1,c_1,c_2$ from $K_2$.

  For the number of singular configurations, note that we have second
  order equations for $s_7$, hence at most 2 values for the tuple
  $(s_3,c_3,s_7,c_7)$, and $s_2$.  So in general if there are two
  separate roots both for $s_7$ and $s_2$, we get four different
  singularities.
\end{proof}
Similar results can be presented for any $T_i$ but we will not
catalogue them here.

\subsection{Subsystem 367}
\label{sec:subsystem367}

Comparing to examples in \cite{ar01:RCS} it was perhaps intuitively
clear that subsystem $J_{4567}$ produces singularities.  It is a bit
more surprising that there is another subsystem producing
singularities: the one formed by the nodes 3, 6, and 7.

Define
\begin{align*}
  h_1 &:= -p_5+p_1
  = a_6 \big(c_6 c_7-s_6 s_7 \big)+a_7 s_7-a_3 s_3+w_1-b_1 \\
  h_2 &:= -p_6+p_2
  = a_6 \big(s_6 c_7+c_6 s_7 \big)-a_7 c_7+a_3 c_3+w_2-b_2 \\
  h_3 &:= p_{9}  = c_3^2+s_3^2-1 \\
  h_4 &:= p_{12} = c_6^2+s_6^2-1 \\
  h_5 &:= p_{13} = c_7^2+s_7^2-1.
\end{align*}
It is important to note that $h_1,h_2$ contain only angles 3,6, and 7,
therefore only $p_9,\,p_{12},\,p_{13}$ are relevant to them. As parameters we now have not only the lengths $a_3,a_6,a_7$, but also
$b_1,\dots,w_2$ i.e.  the positions of the fixed nodes $A$ and $B$ in
Figure \ref{fig:lengths}. Let
$J_{367}$ be the ideal generated by $h_1,\dots,h_5$. We will proceed in a similar way as with the subsystem $J_{4567}$.

First we will consider the singularities of the subsystem $J_{367}$
using the following product order:
\begin{equation}
  \label{eq:S2}
    J_{367}\cup F_{J_{367}} \subset \Q[(c_3,s_3,c_6,s_6,c_7,s_7),(a_3,a_6,a_7,b_1,b_2,w_1,w_2)]
\end{equation}
The relevant Gr\"obner basis $G$ contains 96 generators of which two are 
especially interesting:
\begin{equation}
\begin{aligned}
  g_{12} &= c_6a_6a_7  \\
  g_{1} &= \prod_{i=1}^{4}z_i \qquad\text{ where}  \\
    & z_1 = (a_3-a_6+a_7)^2-|b-w|^2 \\
    & z_2 = (a_3+a_6+a_7)^2-|b-w|^2 \\
    & z_3 = (a_3+a_6-a_7)^2-|b-w|^2 \\
    & z_4 = (a_3-a_6-a_7)^2-|b-w|^2 .
\end{aligned}
    \label{eq:U_z}
\end{equation}
The latter one gives us the singular variety $S_{J_{367}}$.
\begin{thm}
\label{thm:z_i}
 The singular variety of $J_{367}$ is
\[
      S_{J_{367}}=\mathsf{V}(\langle g_1\rangle).
\]
\end{thm}
\begin{rem}
  It is worth noting that, contrary to the linear constraints $t_i$ in
  Theorem \ref{thm:t_i} related to $J_{4567}$, the $z_i$ in Theorem
  \ref{thm:z_i} give {\em quadratic} constraints $z_i=0$ related to $J_{367}$ and
  have the interpretation ``$|a_3\pm a_6 \pm a_7| =$ distance between
  the fixed points A and B''. Furthermore, again the factors $z_i$ give the irreducible decomposition of the singular variety.
\end{rem}
Since $a_i>0$, we get $c_6=0$ from $g_{12}=0$. This simplifies computations considerably. Let us define
\[
U:=\langle J_{367},\, F_{J_{367}},\, c_6 \rangle.
\]
The prime decomposition of $U$ turns out to have 8 components:
\begin{equation*}
  \sqrt{U} = U_1 \cap \dots \cap U_8.
\end{equation*}
Inspecting the generators of each of $U_i$, it is noticed that the
ideals $U_k,\quad k=5\dots 8$ contain generators which imply $a_i=0$
for some $i$.  Hence those are discarded as non-physical and we are left with
4 ideals:
\begin{align*}
  U_{1} &= \langle u_1 ,\,u_2,\, c_7^2+s_7^2-1,\, c_6,\, s_6-1,\, s_3+s_7,\, c_3+c_7 \rangle \\
  U_{2} &= \langle u_1 ,\, u_2 ,\, c_7^2+s_7^2-1,\, c_6,\, s_6+1,\, s_3+s_7,\, c_3+c_7 \rangle \\
  U_{3} &= \langle u_1 ,\,u_2 ,\, c_7^2+s_7^2-1,\, c_6,\, s_6+1,\, s_3-s_7,\, c_3-c_7 \rangle \\
  U_{4} &= \langle u_1 ,\, u_2 ,\, c_7^2+s_7^2-1,\, c_6,\, s_6-1,\, s_3-s_7,\, c_3-c_7 \rangle\\[2mm]
   \textrm{where}&\quad
   \begin{cases}
        u_1=-s_6c_7a_6-c_3a_3+c_7a_7+b_2-w_2\\
        u_2=s_6s_7a_6+s_3a_3-s_7a_7+b_1-w_1.
    \end{cases}
\end{align*}
With these, we continue studying the whole system $J\cup F_J$.  Each $U_i$
will lead to a different case with $s_6=\pm 1$, $s_3=\pm s_7$,
$c_3=\pm c_7$. Let us look for example the ideal $U_1$.\footnote{As
with $J_{4567}$ and $T_5$, the other cases are completely similar and
we will comment them shortly.}  This gives
\begin{equation}
\label{S2_U1}
\begin{aligned}
  s_6 &= 1, \\
  c_7 &= \frac{b_2-w_2}{a_6-a_3-a_7},  \\
  s_7 &= \frac{b_1-w_1}{a_3-a_6+a_7},  \\
  c_3 &= -c_7,\\
  s_3 &= -s_7.
\end{aligned}
\end{equation}
We should expect to run into an equation $z_i=0$ for some $i$, where
the expressions $z_i$ are given in \eqref{eq:U_z}.  Combined with
$c_7^2+s_7^2-1=0$ the equations \eqref{S2_U1} give $z_1=0$.  Likewise,
$U_i$ implies $z_i=0$ for $i=2,3,4$.
\begin{rem}\label{rem:U_i_physical}
  The condition $z_2=0$ is physically a redundant case: it means that
  the system can barely reach from $A$ to $B$ when the subsystem of the rods
  $a_3,a_6,a_7$ is fully stretched, i.e. it has no room to move.
  Therefore also $U_2$ corresponds to a rather trivial case.  See also
  Remark \ref{rem:T_i_physical}.
\end{rem}
Using  $U_1$ we can now eliminate the variables corresponding to angles 3, 6, and 7.  Doing the substitutions in $J\cup F_J$ we are left with the following generators.
\begin{equation}
\label{lideaali}
\begin{aligned}
  L=&\ \langle L_1\cup L_2\rangle,\\
   L_1\quad: &\quad
  \begin{cases}
    l_1= a_2(-c_1c_2+s_1s_2)+c_1a_1+s_7a_3-b_1\\
    l_2= a_2(-s_1c_2-c_1s_2)+s_1a_1-c_7a_3-b_2\\
    l_3= c_1^2+s_1^2-1\\
    l_4= c_2^2+s_2^2-1,
  \end{cases} \\
   L_2\quad : &\quad
  \begin{cases}
    l_5= a_4(s_4c_5+c_4s_5)+c_5a_5+s_7(a_6-a_7)\\
    l_6= a_4(c_4c_5-s_4s_5)-s_5a_5+c_7(a_6-a_7)\\
    l_7= c_4^2+s_4^2-1\\
    l_8= c_5^2+s_5^2-1,
  \end{cases}
\end{aligned}
\end{equation}
where the $s_7,\, c_7$ are no longer variables, but known expressions
from \eqref{S2_U1} and kept here only for clarity of notation.
\begin{rem}
  Before working on $L_1$ and $L_2$ we comment briefly on the other
  $U_i$ cases.  Introduce $L_3$ and $L_4$:
  \begin{align*}
     L_3: &
    \begin{cases}
      & a_2(-c_1c_2+s_1s_2)+c_1a_1-s_7a_3-b_1=0\\
      & a_2(-s_1c_2-c_1s_2)+s_1a_1+c_7a_3-b_2=0\\
      & c_1^2+s_1^2-1=0\\
      & c_2^2+s_2^2-1=0
    \end{cases} \\
     L_4: &
    \begin{cases}
      & a_4(s_4c_5+c_4s_5)+c_5a_5-s_7(a_6+a_7)=0\\
      & a_4(c_4c_5-s_4s_5)-s_5a_5-c_7(a_6+c_7)=0\\
      & c_4^2+s_4^2-1=0\\
      & c_5^2+s_5^2-1=0.
    \end{cases}
  \end{align*}
  Had we used $U_2$ instead of $U_1$, we would end up with the system
  $L_1,\,L_4$.  Likewise, $U_3$ would give the system $L_3,\,L_2$, and
  $U_4$ would give the system $L_3,\,L_4$.  Yet another point of view
  is, that $s_6=\pm 1$ picks between $L_2$ and $L_4$, while
  $(c_3,s_3)=\pm(c_7,s_7)$ picks between $L_1$ and $L_3$.  More
  precisely, $s_6=1$ ($s_6=-1$) gives $L_2$ ($L_4$), and
  $(c_3,s_3)=(-c_7,-s_7)$ gives $L_1$.  The choice
  $(c_3,s_3)=(c_7,s_7)$ would give $L_3$.
\end{rem}
Continuing with $L_1$ and $L_2$, we notice that $L_2$ contains only
the variables $c_5,s_5,c_4,s_4$ (angles 4 and 5), has 4 equations and
4 variables hence is expected to have a finite solution set and will
be handled analogously to the ideal $K_2$ in \eqref{kideaali}.
Calculating its Gr\"obner basis $G$ in the ring
\[
\Q(a_4,a_5,a_6,a_7)[(c_4,c_5,s_5,c_7,s_7),(s_4)]
\]
we obtain 12 generators, the first one being
\[
  g_1 = 2a_4a_5s_4 + a_4^2+a_5^2-a_6^2+2a_6a_7-a_7^2.
\]
Hence $s_4$ can be explicitly solved:
\begin{equation}
\label{S2_U1_s4}
  s_4 = \frac{a_4^2+a_5^2-a_6^2+2a_6a_7-a_7^2}{-2a_4a_5}.
\end{equation}
The other generators are too messy to be of much use. Then using the formula $c_4^2=1-s_4^2$ we get

\begin{align}\label{S2_U1_c4}
  c_4^2 &= -\frac{(a_4+a_5-a_6+a_7)(a_4-a_5+a_6-a_7)(a_4-a_5-a_6+a_7)(a_4+a_5+a_6-a_7)}{4a_4^2a_5^2} \notag \\
  &= -\frac{t_7 t_5 t_6 t_8}{4a_4^2a_5^2}.
\end{align}
The product term in the numerator has to be nonpositive, in order to
have any real solutions:
\begin{equation}
  t_5 t_6 t_7 t_8 \le 0.
\end{equation}
After solving $s_4,c_4$ we can proceed to solve $s_5$ and $c_{5}$.
For this we use the ordering
\[
\Q(a_4,a_5,a_6,a_7)[c_5,s_5,c_4,s_4,c_7,s_7]
\]
and pick the two relevant equations from the corresponding Gr\"obner basis:
\begin{align*}
  (-a_6+a_7)s_5 -a_4c_4s_7 +a_4s_4c_7 +a_5c_7 &=0 \\
  (-a_6+a_7)c_5 -a_4c_4c_7 -a_4s_4s_7 -a_5s_7 &=0,
\end{align*}
which are linear equations for $s_5,c_5$, provided $a_6\neq a_7$.
\begin{rem}
  In the case $a_6=a_7$ the situation is different: $L_2$ then
  decomposes into 3 prime ideals, of which only one is physically
  feasible and gives a singularity only if $a_4=a_5$.  Thence this is
  a rather special case and will not be considered further here.
\end{rem}
The subsystem $L_2$ is now fully solved.  
Moving on to $L_1$, we will
see that the analysis is very similar to that of $K_1$ from
\eqref{kideaali}.  Therefore we will skip some details.
After forming the Gr\"obner basis of $L_1$ in the ring
\[
\Q(b_1,b_2,a_1,a_2,a_3,c_7,s_7)[c_1,s_1,c_2,s_2]
\]
with respect to the lexicographic ordering, we get 
for $s_2$, after simplifications, the relation
\begin{align}
  s_2^2 &= \frac{n_{3}(4a_1a_2-n_3)}{4a_1^2a_2^2},\\
  &\qquad\text{ where} \quad
    n_3 = |b|^2+2a_3(b_2c_7-b_1s_7)-(a_1-a_2)^2+a_3^2 \notag
\end{align}
Again for the real solutions the numerator has to be
nonnegative
\begin{equation}
n_{3}(4a_1a_2-n_3)\ge 0
\end{equation}
We can now solve $c_{2}$, $s_{1}$ and $c_1$, provided
their coefficients are nonzero, from the linear equations
\begin{align*}
  &2a_1a_2n_4c_2-4a_1^2a_2^2s_2^2+r_1=0, \\
  &-2a_1n_4s_1+r_2+r_3=0, \\
  &({a_1^2-a_2^2})c_1+r_4=0.
\end{align*}
where 
\[
   n_4=|b|^2+a_3^2+2a_3(b_2c_7-b_1s_7)
\]
and $r_i$ are lengthy, yet polynomial, expressions in the
parameters, apart from $r_4$ which depends on
$s_1,s_2,c_2$ as well.  (See the appendix.)

What about the cases $n_4=0$ and/or $a_1=a_2$?
It can be shown, as with $d_0$, that the condition $n_{3}(4a_1a_2-n_3)\ge 0$ is equivalent to
\[
 (a_1-a_2)^2 \le n_4 \le (a_1+a_2)^2.
\]
Therefore, if $a_1\neq a_2$ then $n_4\neq 0$ and the equations above are
sufficient.  The case $a_1=a_2$, $n_4\neq 0$ does not essentially change the
situation: we still have a quadratic equation for $s_2$, and linear ones 
for the others, with a different coefficient for $c_1$.

The remaining case  $a_1=a_2$, $n_4=0$ is analogous to the $n_2=0$ case
within $J_{4567}$ and likewise will not be pursued further.

\begin{thm}
  \label{thm:2}
  Let us suppose that the parameters $a$,$b$,$w$ satisfy the following
  conditions:   
  \begin{align}
    & a_6\neq a_7 \notag \\
    & n_4 \neq 0 \notag \\
    & n_{3}(4a_1a_2-n_3)\geq 0 \label{U1_cond_1} \\
    & t_7 t_5 t_6 t_8\leq 0 \label{U1_cond_2}.
  \end{align}
  Then $V_{a,b}$ contains at least 2 singular points. If the
  inequalities are strict we get in general at least 4 singular
  points.
\end{thm}

Similar results can be represented for any $\mathsf{V}(U_i)$ but we will not
catalogue them here.

\begin{proof}
  The last two conditions are due to univariate second degree
  polynomial equations, which have real solutions if and only if
  \eqref{U1_cond_1} (for $s_2$) and \eqref{U1_cond_2} (for $c_4$) are
  fulfilled.  The first condition is needed for the other variables to
  be determined uniquely: $s_3,c_3,s_6,c_6,s_7,c_7$ from
  $\mathsf{V}(U_1)$, $s_4,s_5,c_5$ from $L_2$, and $s_1,c_1,c_2$ from
  $L_1$.

  For the number of singular configurations, note that we have second
  order equations, hence at most 2 values, for $c_4$ and $s_2$.  So in
  general if there are two separate roots both for $c_4$ and $s_2$, we get
  four different singularities.
\end{proof}

\subsection{Two special cases with symmetry}
Let us look more closely at two special cases: $a_4=a_6,\,a_5=a_7$, 
and either $a_4=a_5$ or $a_4\neq a_5$.

\subsubsection{The case $a_4\neq a_5$}
\label{sec:a4_not_a5}

Motivated by the original benchmark values \cite{sc90:MSH} we give the
following
\begin{lem}\label{lem:1}
  When $a_4=a_6$ and $a_5=a_7$, there is a relation between the angles
  4 and 6: either $y_6=-y_4$ or $y_6=y_4+\pi$.  Furthermore, if also
  $a_4\neq a_5$, the angle $y_7$ variables, i.e. $c_7,s_7$, are
  uniquely determined from $c_4,s_4,c_5,s_5$.
\end{lem}
\begin{proof}
  Looking for relations between solely angles 4 and 6, we substitute
  $a_4=a_6$ and $a_5=a_7$ to the subsystem $J_{4567}$ and formulate a
  suitable elimination ideal.  In ideal-theoretic language, we define
  \begin{align*}
    & r_1:=a_4\big(s_4c_5+c_4s_5\big)+a_5c_5
    -a_4\big(c_6c_7-s_6s_7\big)-a_5s_7\\
    & r_2:=a_4\big(c_4c_5-s_4s_5\big)-a_5s_5
    +a_4\big(s_6c_7+c_6s_7\big)-a_5c_7\\
    & r_{i+2}=c_{i+3}^2+s_{i+3}^2-1,\quad i=1,\dots,4,
  \end{align*}
  where $r_i=q_i$ with substitutions $a_4=a_6$ and $a_5=a_7$, and
  investigate the ideal $I:=\langle r_1,\ldots,r_6\rangle$ in the ring
  \[
  \Q(a_4,a_5,a_6,a_7)[(c_5,s_5,c_7,s_7),(c_4,s_4,c_6,s_6)].
  \]
  Calculating the elimination ideal
  $I_{4,6}:=I\cap\Q[c_4,s_4,c_6,s_6]$ we get
  \begin{eqnarray*}
    I_{4,6}= \langle s_4+s_6,c_6^2+s_6^2-1,c_4^2+s_4^2-1\rangle.
  \end{eqnarray*}
  Calculating the prime decomposition of $\sqrt{I_{4,6}}$ we get
  \begin{equation*}
    \sqrt{I_{4,6}}=\langle c_6^2+s_6^2-1,c_4-c_6,s_4+s_6\rangle
    \cap\langle c_6^2+s_6^2-1,c_4+c_6,s_4+s_6\rangle.
  \end{equation*}
  Since $I_{4,6}\subset I \subset J \subset J\cup F_J$, we have
  $$
  \mathsf{V}(I_{4,6}) \supset  \mathsf{V}(J\cup F_J).
  $$
  From these prime ideals we can see that everywhere in $\mathsf{V}(I_{4,6})$,
  and therefore in the variety of the singularities of the whole
  system as well, $s_6=-s_4$ and either $c_6=c_4$ or $c_6=-c_4$.
  These translate into two possible relations between the angles $y_4$
  and $y_6$.
  \begin{equation}
    (c_6,s_6)=(c_4,-s_4) \Leftrightarrow y_6=-y_4, \qquad
    (c_6,s_6)=(-c_4,-s_4) \Leftrightarrow y_6=y_4+\pi.
  \end{equation}
  This proves the first claim.  If we take into account either one of
  the prime ideals of $\sqrt{I_{4,6}}$ in $I$ and calculate the
  Gr\"obner bases we get ideals where $c_7$ and $s_7$ depend linearly
  on $c_4$, $s_4$, $c_5$ and $s_5$, and can be explicitely solved, as
  we will show next to prove the latter claim of the lemma.  For the
  case $(s_6,c_6)=(-s_4,-c_4)$ we get
  \begin{equation}
    \begin{cases}
      c_7 &= -s_5\\
      s_7 &= c_5
    \end{cases} \text{ which imply }
    y_7=y_5+\frac{\pi}{2}.
  \end{equation}
  For the case $(s_6,c_6)=(-s_4,c_4)$ the expressions are, albeit linear,
  slightly more complicated:
  \begin{align*}
    &c_7\big(a_4^2(s_4^2-c_4^2)-a_5(2a_4s_4+a_5)\big)+s_7\big(2a_4(a_5+a_4c_4s_4)\big)-s_5\big((a_4^2+a_5^2)-2a_4a_5s_4\big) =0 \\
    -&c_7\big(2a_4^2c_4s_4\big)+s_7\big(a_4^2(c_4^2-s_4^2)+a_5^2)\big)+(a_4^2-a_5^2)c_5-2a_4a_5s_5c_4 =0.
  \end{align*}
  We prove that these indeed determine $c_7,s_7$: all we need to do is
  check that the determinant of the coefficient matrix $A$ of the
  linear equations does not equal zero:
  \begin{equation*}
    A:=
    \begin{pmatrix}
      a_4^2(s_4^2-c_4^2)-a_5(2a_4s_4+a_5) & 2a_4(a_5+a_4c_4s_4) \\
      -2a_4^2c_4s_4 & a_4^2(c_4^2-s_4^2)+a_5^2
    \end{pmatrix},\text{ prove }
    \det(A) \neq 0.
  \end{equation*}
  Now det$(A)$ simplifies due to $c_4^2+s_4^2=1$, resulting in
  \begin{align*}
    \det(A)
    &=2a_4a_5(a_4+a_5)(a_4-a_5)s_4+(a_4-a_5)(a_4+a_5)(a_4^2+a_5^2)
  \end{align*}
  Let us then consider det$(A)$ as a function of $s_4$.  Since
  $s_4\in[-1,1]$, $\det(A):[-1,1]\mapsto\R$.  Clearly if $a_4=a_5$,
  $\det(A)\equiv 0$ so we need to assume $a_4\neq a_5$.  Set
  $$
  h(s_4):=\frac{\det(A)}{(a_4+a_5)(a_4-a_5)}=2a_4a_5s_4+(a_4^2+a_5^2)
  $$
  and inspect when $h=0$.  Since $a_4>0$ and $a_5>0$ the linear function
  $h$ has its minimum at $-1$.
  $$
  h(-1)=a_4^2+a_5^2-2a_4a_5=(a_5-a_4)^2>0.
  $$
  This proves $h\neq 0$ always, therefore under the assumption $a_4\neq
  a_5$ also $\det(A)\neq 0$ as claimed.
\end{proof}

\subsubsection{The case $a_4=a_5$}
\label{sec:a4_is_a5}

We study the special case $a_4=a_5=a_6=a_7$, whence the
  4567-subsystem is capable of ``buckling'' in more complicated ways,
  thereby producing further interesting configurations.  This resembles then
  the net example in \cite{ar01:RCS}.

\begin{figure}[htb]
  \centering
     \epsfig{file=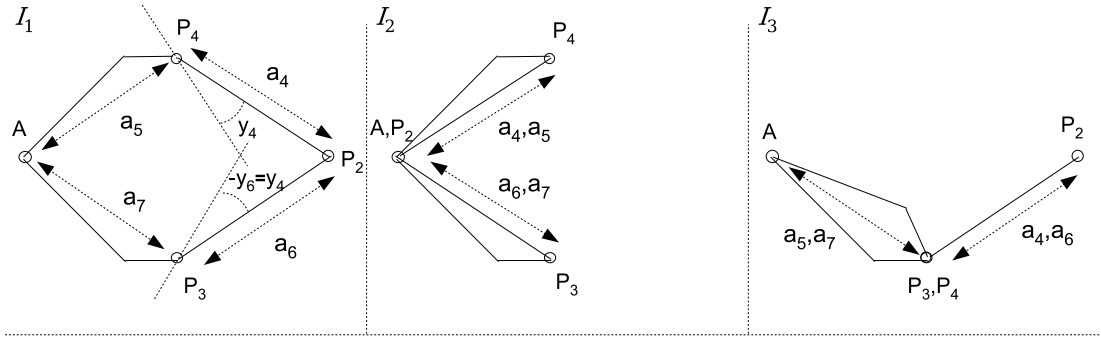, width=0.9\textwidth, height=0.3\textwidth}
  \caption{The configurations corresponding to $I_1,\,I_2,\,I_3$ in
    the case $a_4=a_5=a_6=a_7$.}
  \label{fig:a4567_equal}
\end{figure}

Let us see how $J_{4567}$ simplifies with substitutions $a_4=a_5=a_6=a_7$.
Note that the assumptions of Lemma \ref{lem:1} considering $y_7$ are
no longer valid.  Let
$$
I := J_{4567} \text{ with }a_4=a_5=a_6=a_7 \text{ and }s_6=-s_4
$$
and compute its prime decomposition.  This results in
\begin{multline}\label{eq:a4_is_a5}
  \sqrt{I} = I_1 \cap I_2 \cap I_3 \quad \text{ with generators } \\[3mm]
  I_1 =
  \begin{cases}
    s_4^2+c_6^2-1,\\
    c_4-c_6,\\
    c_7^2+s_7^2-1,\\
    s_5+c_7s_4-s_7c_6,\\
    c_5-c_7c_6-s_7s_4
  \end{cases} \qquad
  I_2 =
  \begin{cases}
    c_6,\\
    s_4+1,\\
    c_4,\\
    c_7^2+s_7^2-1,\\
    c_5^2+s_5^2-1
  \end{cases} \qquad
  I_3 =
  \begin{cases}
    s_4^2+c_6^2-1,\\
    c_4+c_6,\\
    c_7^2+s_7^2-1,\\
    s_5+c_7,\\
    c_5-s_7
  \end{cases}
\end{multline}
Each of these has a geometrical interpretation, see Figure
\ref{fig:a4567_equal}.  $I_2$ corresponds to $y_4=-\pi/2, y_6=\pi/2$
which means that nodes $A$ and $P_2$ coincide.  This is like the
$T_5$ situation.  Indeed, the ideal $J\cup F_J\cup I_2$ turns out to
be exactly $T_5$ with the extra condition $a_4=a_5$.  Although it is
not immediately apparent but in that situation there also arises a new
condition among the parameters: $a_3=|b-w|$, i.e. ``$a_3$ equals the distance between $A$ and $B$''.
Note that here the Fitting ideal $F_{J_{4567}}$ has not been used at
all, contrary to the $T_5$ calculations.

$I_3$ corresponds to $y_6=y_4+\pi$ and $y_5=y_7-\pi/2$ so that now nodes $P_3$ and $P_4$ coincide.  Then again, $I_1$
corresponds to $y_6=-y_4$ and $y_5=y_6+y_7$, which interestingly is
{\em not a singularity} but merely expressing a symmetry in the system
due to $a_4=a_5=a_6=a_7$.

\subsection{Other subsystems}

Now contemplating Figure \ref{fig:lengths} we see that it would be possible to find other singularities by analysing still other subsystems. For example the subsystem  corresponding to rods 3, 4 and 5 is by symmetry similar to subsystem 367: we simply exchange the roles of variables and parameters associated to rods 4 and 6, and 5 and 7. Further we could consider other subsystems formed from
different ``paths'' between the nodes $A,B,O$: i.e. subsystems
$J_{123},J_{1245},J_{1267}$. Again by symmetry the system $J_{1267}$ is completely
similar to $J_{1245}$, but cases $J_{123}$ and $J_{1245}$ give new singularities. We checked that in these cases the singular variety is not empty, and that at least for some parameter values we get singular points. 

We did not analyse these cases in detail because computations are quite similar to those given above for subsystems $J_{4567}$ and $J_{367}$.  Hence we did not feel including these would
give significant additional value and therefore left them out to avoid
expanding this quite a long presentation further.

\section{Numerical examples}
\label{sec:examples}

In this section we will calculate numerical examples for both types of
singularities.  Interestingly, the explicit expressions within
$\tilde{G},\,\hat{G}$, as well as in the Gr\"obner bases of $L_1$ and
$L_2$, are unstable for numerical computations.  It is better to use
the original defining equations of $K_1,K_2,L_1,L_2$ in the
computations.  We shall not explore this stability issue here due to
its non-relevance for the present context.

We present 4 examples:
\begin{enumerate}
\item The original benchmark parameter values, see \cite{testset}.  We
  show that then the system is avoiding
  singularities.\footnote{Thereby validating its benchmark status.
    That is, the numerical difficulties encountered there are indeed
    due to the ``numerical stiffness'' of the problem, not to a nearby
    singularity.}
\item We explore how should $a_1,a_2$ be changed in order to have
  $J_{4567}$ type singularities in the system.  Here we have an
  interpretation for the result: the lengths $a_1,a_2$ must be such
  that the subsystem 4567 has room for a certain kind of ``buckled''
  configuration.
\item We explore how should $b_1,a_1,a_2$ be changed in order to have
  $J_{367}$ type singularities in the system.
\item A special case which shows a rational solution, that is
  $c_i,s_i\in\Q$ for all $i$. This shows unambiguously that we can find singular points because in this case there are no numerical errors related to floating point computations.
\end{enumerate}

\subsection{Original values}
\label{Exam:1}
In this example, we will use the original values for the parameters
$a_i,b_i$ and show that the system then has no
singularities.  The original parameters used in the benchmark tests
\cite{sc90:MSH,ha-wa91:SODE2,testset} are
\begin{align} \label{orig_a_b}
  & a_1=0.007\quad a_2=0.028  \quad a_3=0.035\quad
  a_4=0.020  \quad a_5=0.040\quad a_6=0.020  \quad a_7=0.040 \nonumber \\
  & b_1=-0.03635\quad b_2=0.03273\quad  w_1=-0.06934\quad w_2=-0.00227.
\end{align}
Since $a_7=a_5$ and $a_6=a_4$, we have $t_4=t_6=0$ (and $t_1<0$,
$t_5<0$) so we could have an $J_{4567}$ singularity: $T_3$ or $T_5$.  
\begin{rem}
  Interpretation: both $T_3$ and $T_5$ describe a situation where the
  4567 system has 'collapsed' into a 1-dimensional object.  The ideal
  $K_2$ tells us how $a_3$ restricts the possible attitudes of 4567.  In
  $T_5$ the centre node $P_2$ has been pushed in, in $T_3$ it has been
  pulled out.
\end{rem}
Let us look more closely first at $T_5$, say, and check the conditions
\eqref{S1_T5_D_condition} and \eqref{S1_T5_E_condition}.  The first one
is fulfilled.
For $E$ we first need to solve $c_3,s_3$ from $\mathsf{V}(K_2)$.
Their solutions are
\begin{multline}
  \label{eq:T5c3s3}
  (c_3,s_3,c_7,s_7)\in \{ (0.4299535996, \,-0.9028509856, \,-0.9975812008, \,0.06951077517),\\
  (0.9266735994, \,-0.3758670513, \,-0.1283212011, \, 0.9917326602)\}
\end{multline}
With these $c_3,s_3$ we can compute $E$.  Both sets in
\eqref{eq:T5c3s3} give $E=\mathcal{O}(10^{-5})>0$ and the condition
\eqref{S1_T5_E_condition} is violated, hence there are no
($J_{4567}-$)singularities.  
What about other singularities?  This is answered by the following
\begin{thm}\label{thm:3}
  With the original benchmark parameter values \eqref{orig_a_b}, the
  Andrews' squeezing system has no singularities.
\end{thm}
\begin{proof}
  We now have $a_4=a_6$, $a_5=a_7$ and $a_4\ne a_5$.  Lemma \ref{lem:1} implies
  variables $c_6$, $s_6$, $c_7$, $s_7$, and so $y_6$ and $y_7$ can be
  explicitely solved in terms of $c_4$, $s_4$, $c_5$, and $s_5$.  It
  is then possible to reduce the original system of constraint
  equations, by forgetting the last two equations from \eqref{eq:2},
  and consider
  \begin{eqnarray*}
    \begin{cases}
      a_1\cos(y_1) - a_2\cos(y_1 + y_2) - a_3\sin(y_3) - b_1 &= 0 \\
      a_1\sin(y_1) - a_2\sin(y_1 + y_2) + a_3\cos(y_3) - b_2 &= 0 \\
      a_1\cos(y_1) - a_2\cos(y_1 + y_2) - a_4\sin(y_4 + y_5) - a_5\cos(y_5) - w_1 &= 0 \\
      a_1\sin(y_1) - a_2\sin(y_1 + y_2) + a_4\cos(y_4 + y_5) - a_5\sin(y_5) - w_2 &= 0.
    \end{cases}
  \end{eqnarray*}
  These are equivalent to
  \begin{eqnarray*}
    \begin{cases}
      a_1\cos(y_1) - a_2\cos(y_1 + y_2) - a_3\sin(y_3) - b_1 &= 0 \\
      a_1\sin(y_1) - a_2\sin(y_1 + y_2) + a_3\cos(y_3) - b_2 &= 0 \\
      - a_4\sin(y_4 + y_5) - a_5\cos(y_5)+a_3\sin(y_3)+(b_1 - w_1) &= 0 \\
      a_4\cos(y_4 + y_5) - a_5\sin(y_5) -a_3\cos(y_3)+(b_2- w_2) &= 0 \\
    \end{cases}
  \end{eqnarray*}
  These can be again represented as polynomials.
  \begin{align*}
    & m_1 := a_1  c_1-a_2\big(c_1c_2-s_1s_2\big)-a_3s_3-b_1=0\\
    & m_2 := a_1  s_1-a_2\big(s_1c_2+c_1s_2\big)+a_3c_3-b_2=0\\
    & m_3 := a_1  c_1-a_2\big(c_1c_2-s_1s_2\big)-a_4\big(s_4c_5+c_4s_5\big)
    -a_5c_5-w_1=0\\
    & m_4 := a_1  s_1-a_2\big(s_1c_2+c_1s_2\big)+a_4\big(c_4c_5-s_4s_5\big)
    -a_5s_5-w_2=0\\
    & m_{i+4} := c_i^2+s_i^2-1=0,\quad i=1,\dots,5.
  \end{align*}
  Substituting the original parameter values \eqref{orig_a_b}, as
  rational numbers, into the polynomials $m_i$ we form an ideal
  $I:=\langle m_1,\ldots,m_9\rangle$.  Let $K:=I\cup F_I$, where $F_I$
  is the Fitting ideal of $I$, and inspect $K$ in the ring
  \[
  \Q[(c_1,s_1,c_2,s_2), (c_3,s_3,c_4,s_4,c_5,s_5)].
  \]
  Now it is possible to compute the Gr\"obner basis $G_K$ for $K$
  explicitly (unlike for $J\cup F_J$ in the introduction) and results
  in
  \begin{eqnarray*}
    G_{K}=\langle 1 \rangle.
  \end{eqnarray*}
  This implies $\mathsf{V}(K)=\emptyset$, proving that with these original
  parameter values there are no singularities.
\end{proof}

\subsection{$J_{4567}$ singularity: original values, apart from $a_1,a_2$}
\label{Exam:2}
Let us see how changing $a_1$ and/or $a_2$ might produce $J_{4567}$ type
singularities.  Our analysis reveals that by suitable combinations of
$a_1$ and $a_2$ we can get between zero and four singularities (of
type $J_{4567}$, that is).  The number of singularities is determined by
$c_3,s_3$, and $E$.

\begin{figure}[htb]
  \centerline{
    \epsfig{file=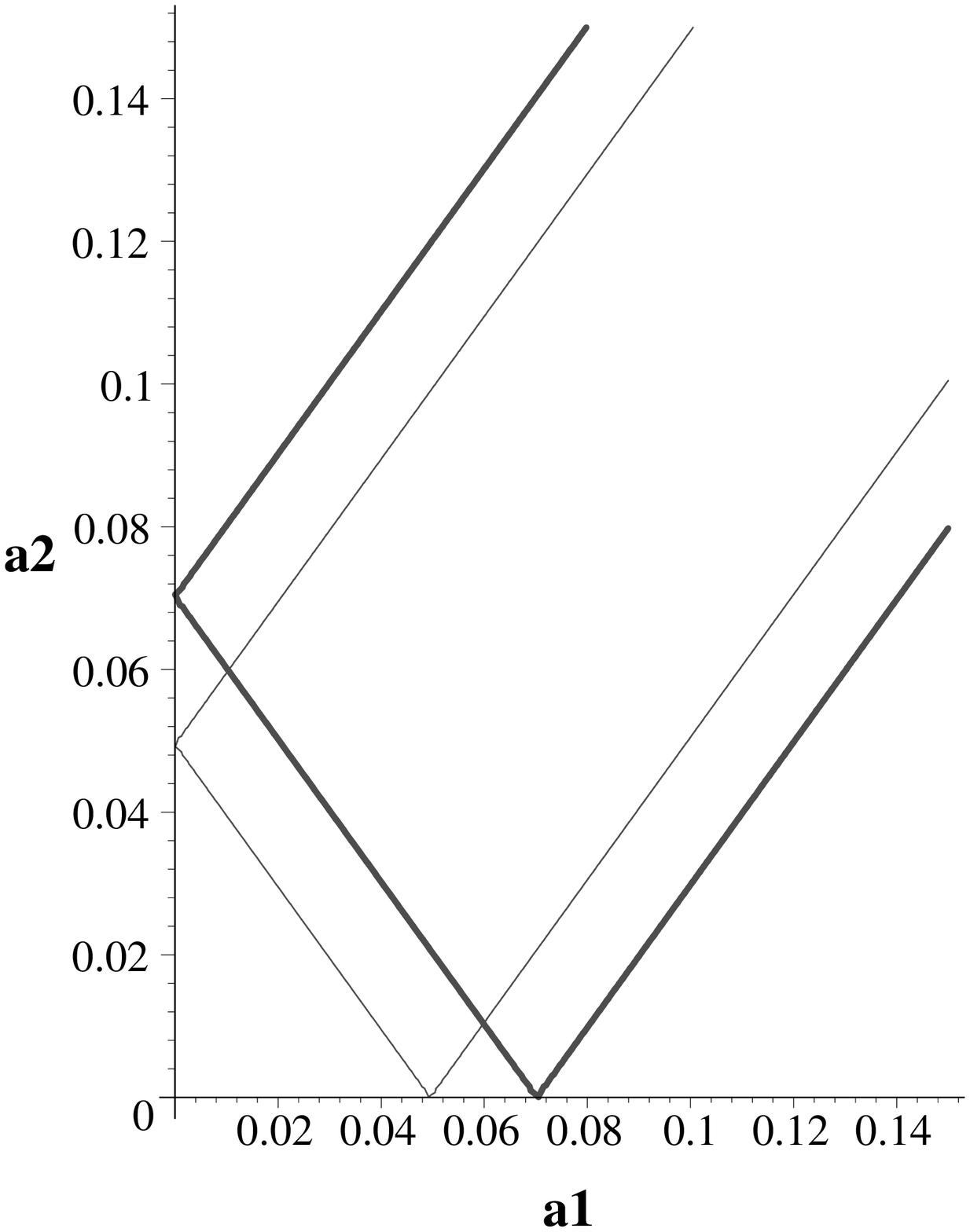,height=7cm,width=7cm} \qquad
    \epsfig{file=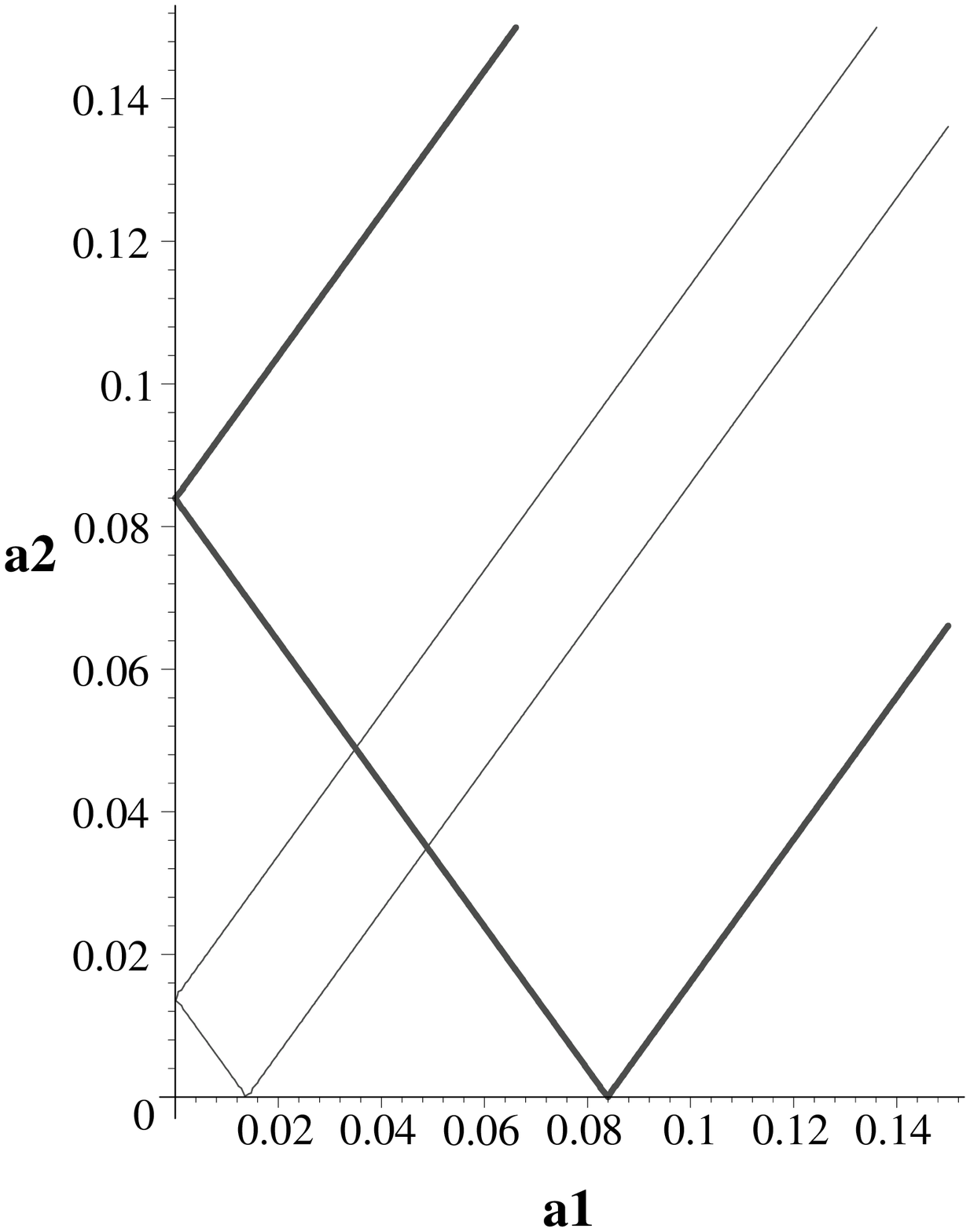,height=7cm,width=7cm}
  }
  \caption{The rectangular lines are $E_1=0$ (thick line) and $E_2=0$
    (thin line), the areas inside each $E_i=0$ line are where $E_i<0$.
    Left panel: $T_5$ case, right panel: $T_3$ case.}
  \label{fig:E_1_and_2}
\end{figure}

Considering $E$ as a function of $a_1,a_2$ we plot the area where
$E\le 0$.  Recall that $E$ depends on $c_3$ as well, and $c_3$ has two
possible values so we get two functions: $E=E_1(a_1,a_2)$ (resp.
$E=E_2(a_1,a_2)$) corresponding to the first (resp. second) value of
$c_3$ from \eqref{eq:T5c3s3}.  See Figure \ref{fig:E_1_and_2} where
the areas inside the rectangular areas are $E_i<0$.
\begin{itemize}
\item no singularities: $E_1>0,E_2>0$.
\item 1 singularity: $E_1=0,E_2=0$, which leads (with $T_5$) to two
  possible values:
  \[
  (a_1=0.05986, \,a_2=0.01035 ), \quad 
  (a_1=0.01035, \,a_2=0.05986 )
  \]
\item 2 singularities: one of $E_1,E_2$ is $< 0$, the other one $>0$.
\item 3 singularities: one of $E_1,E_2$ is $< 0$, the other one $=0$.
\item 4 singularities: $E_1 < 0$, $E_2 < 0$.
\end{itemize}

\begin{figure}[htb]
  \centering
  \begin{tabular}{l|r}
    \epsfig{file=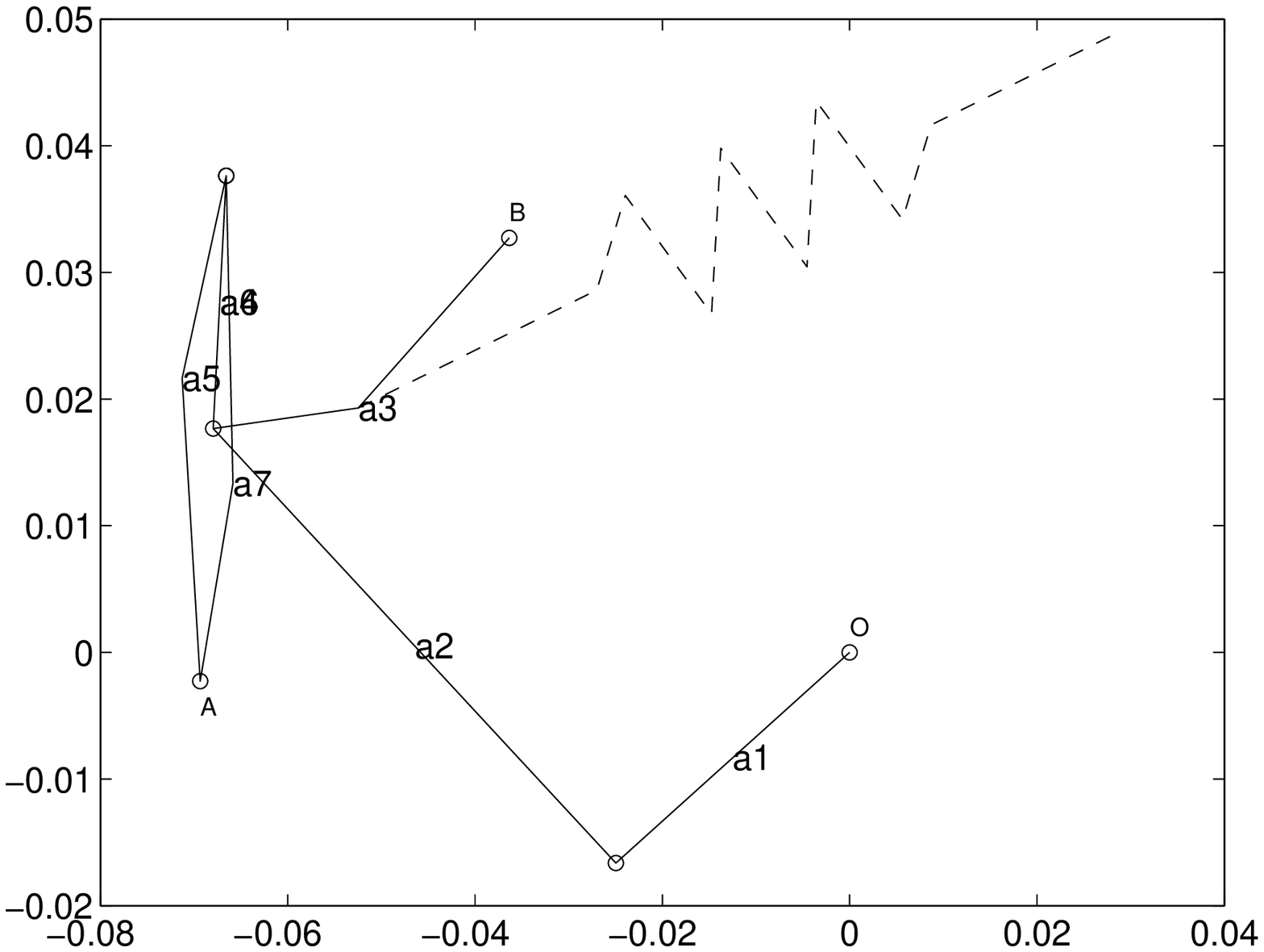, width=0.45\textwidth} &
    \epsfig{file=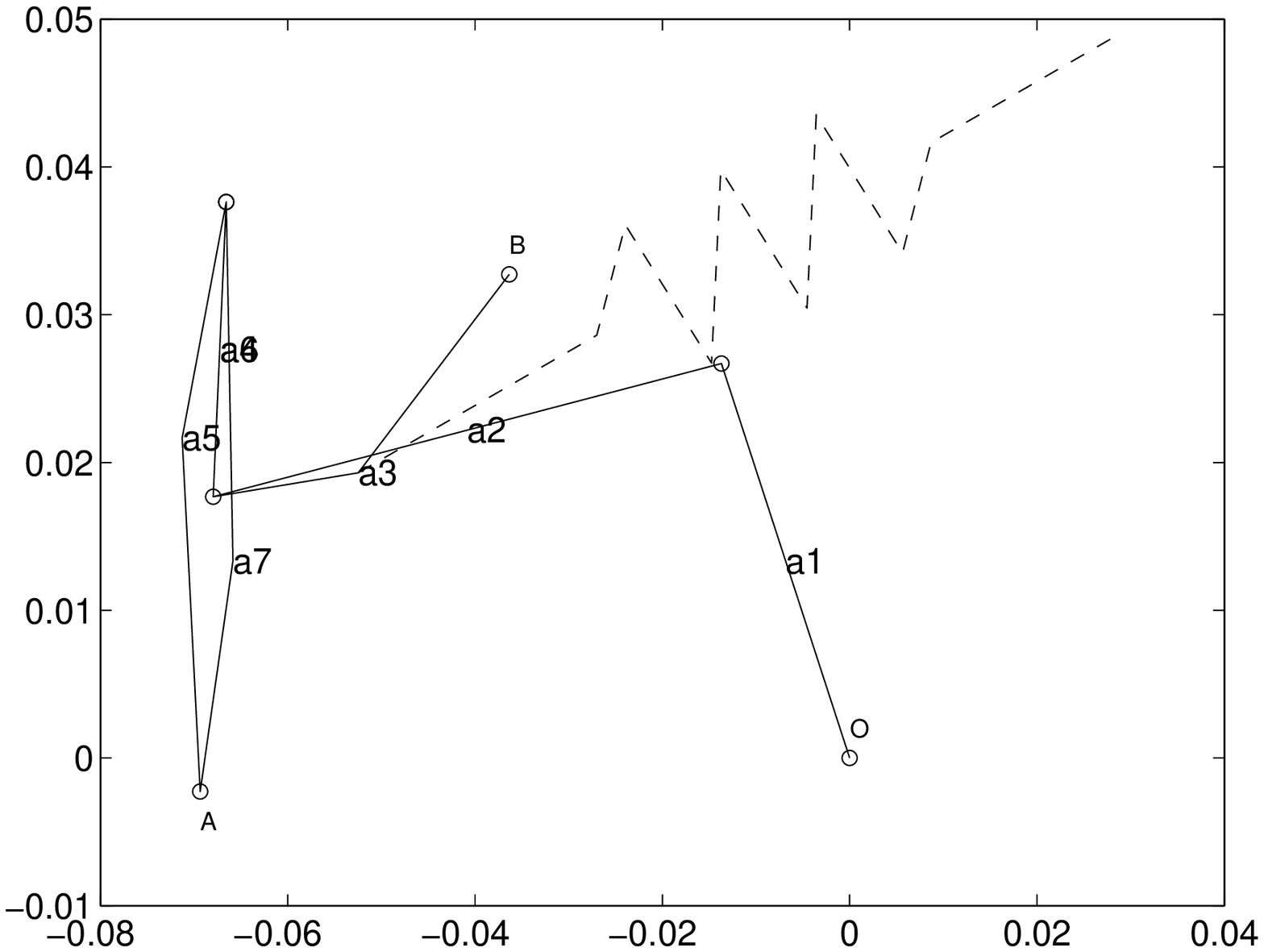, width=0.45\textwidth} \\
    \hline
    \epsfig{file=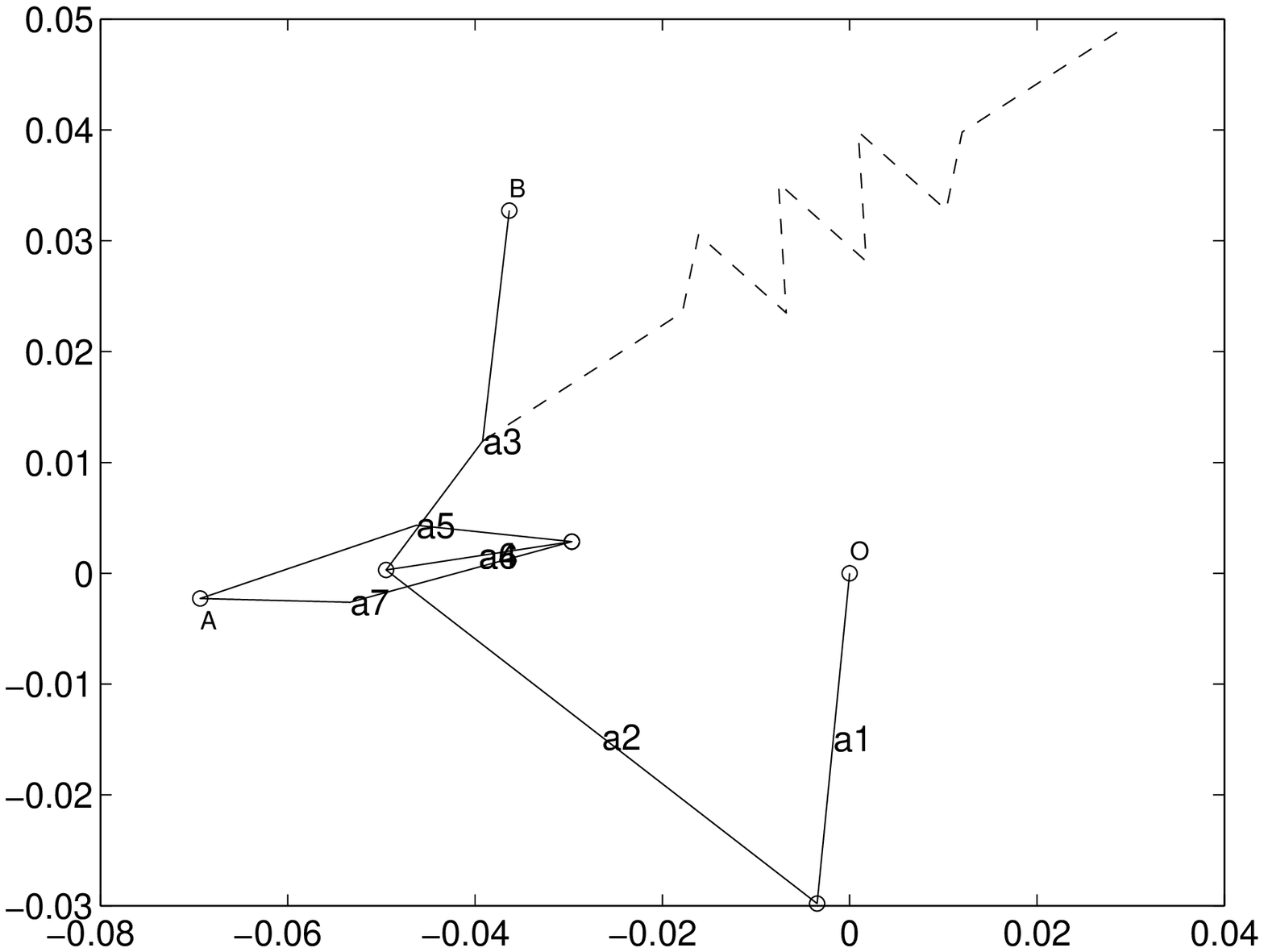, width=0.45\textwidth} &
    \epsfig{file=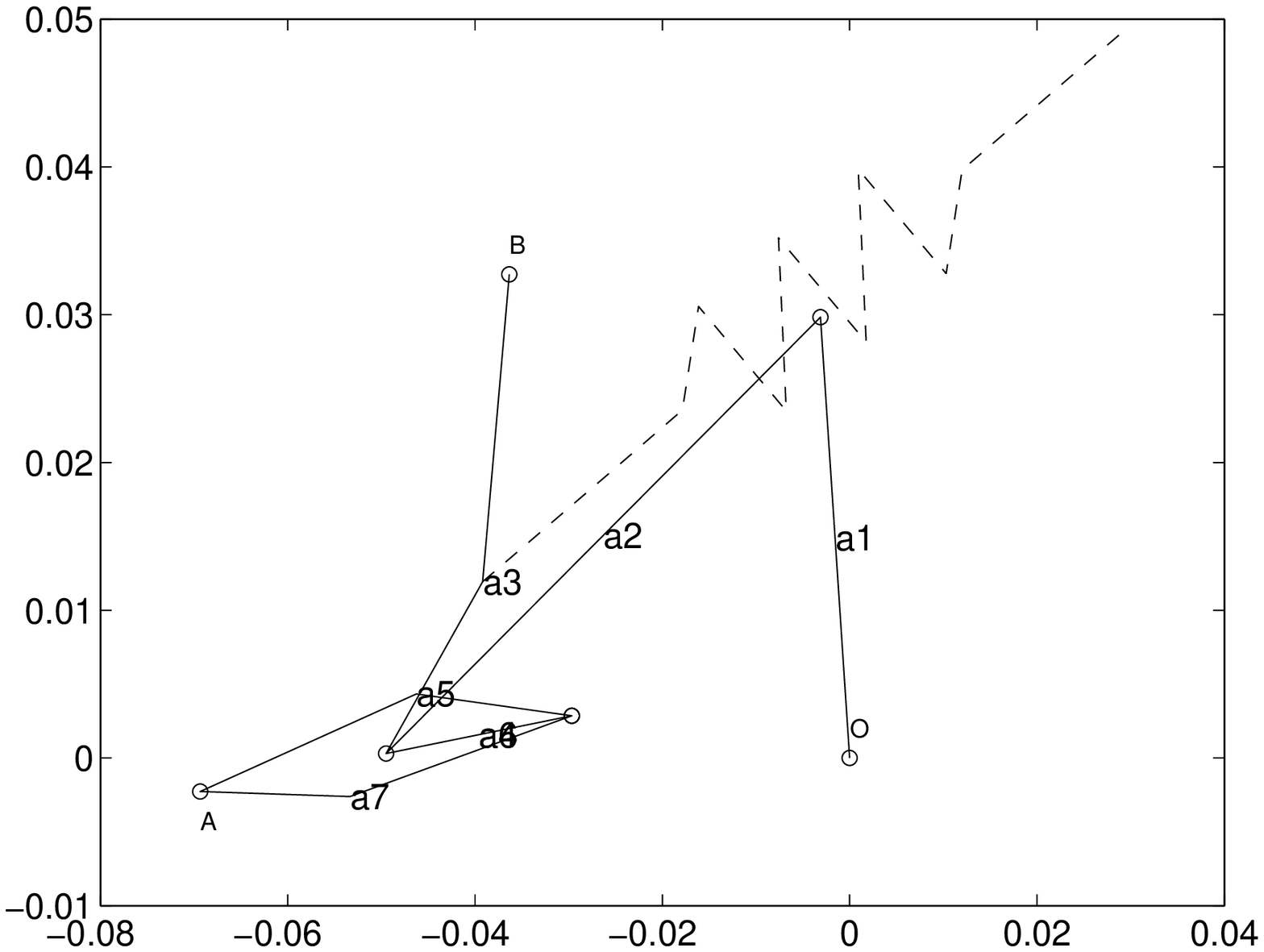, width=0.45\textwidth} \\
    \hline
  \end{tabular}
  \caption{Singular positions (according to $J_{4567}$, $T_5$) when
    $a_1=0.03,\,a_2=0.055$ and $a_3,\dots,a_7$ have the original
    values.  One can see a physical explanation to the singularity:
    the centre node $P_2$ is 'pushed in' so that nodes $P_3$ and $P_4$ coincide.}
  \label{fig:S1_T5_sing}
\end{figure}

\begin{figure}[htb]
  \centering
  \begin{tabular}{l|r}
    \epsfig{file=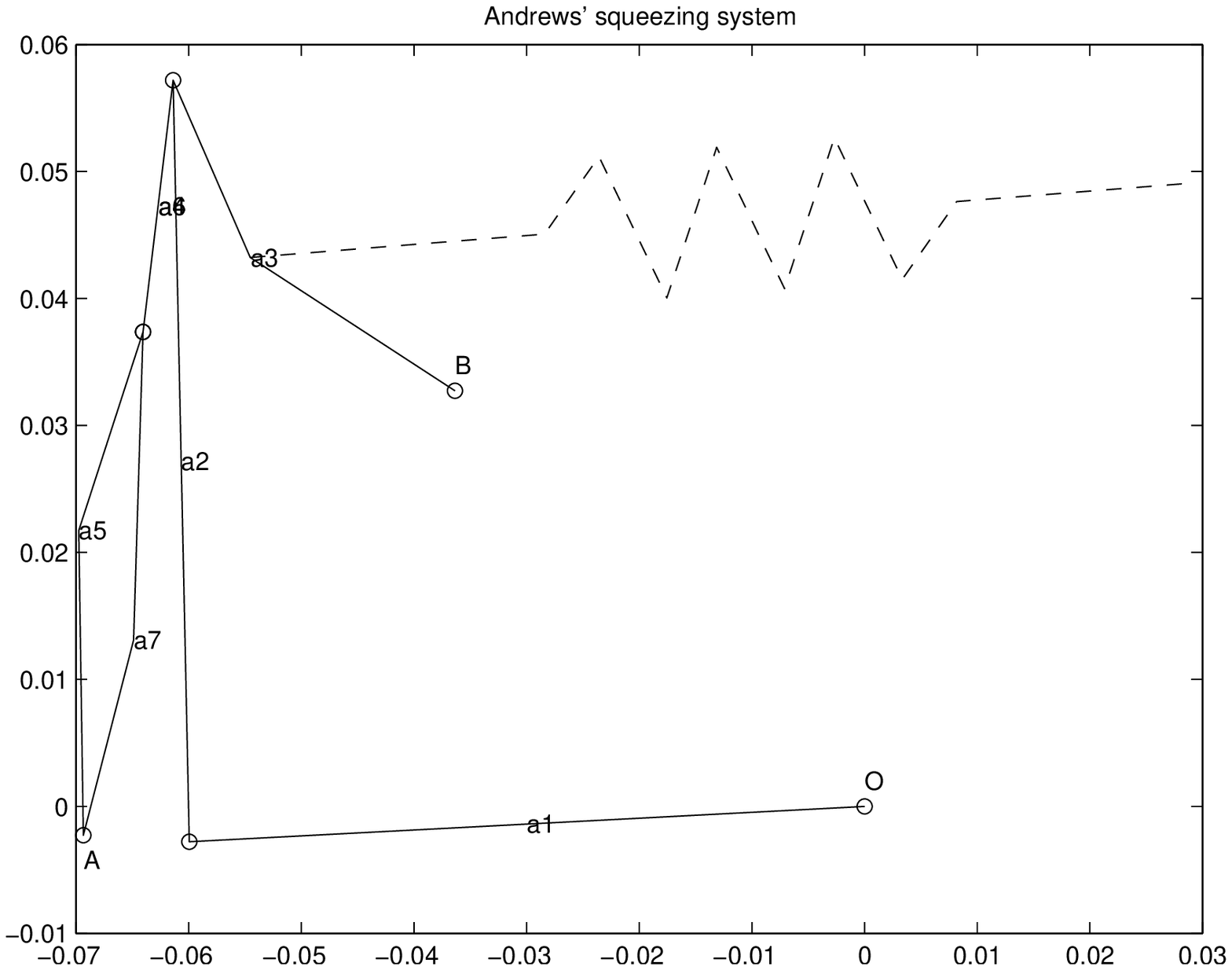, width=0.45\textwidth} &
    \epsfig{file=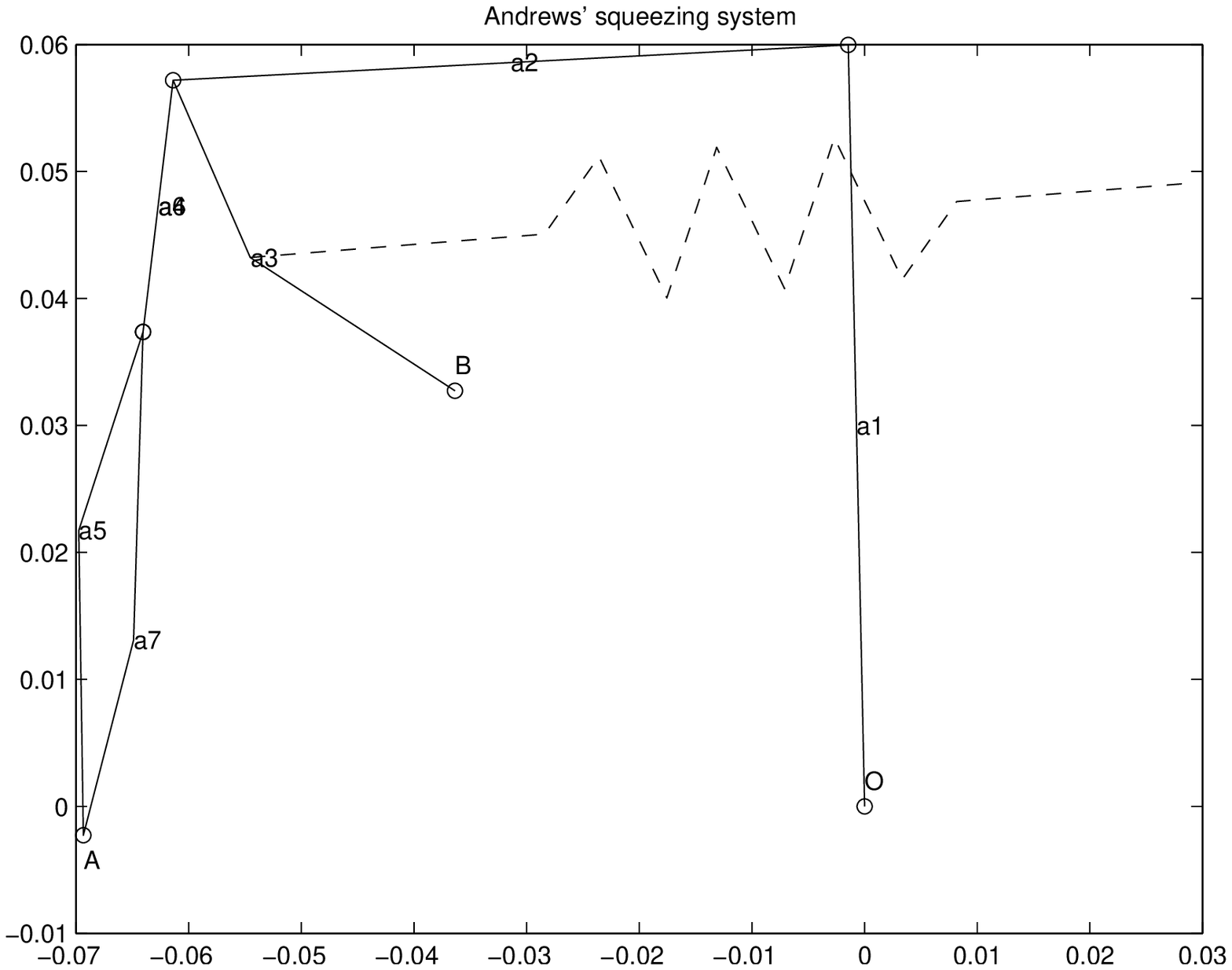, width=0.45\textwidth} \\
    \hline
    \epsfig{file=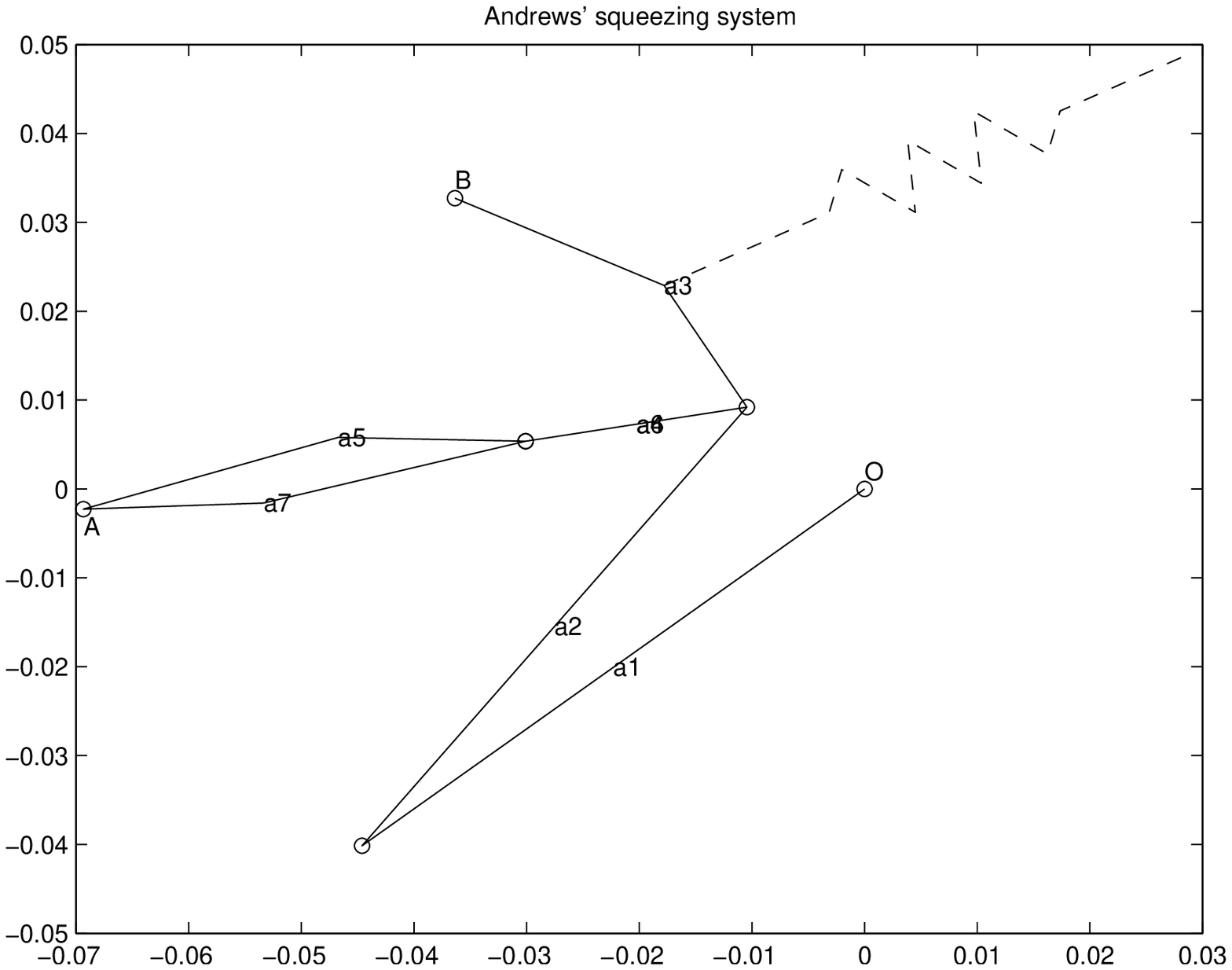, width=0.45\textwidth} &
    \epsfig{file=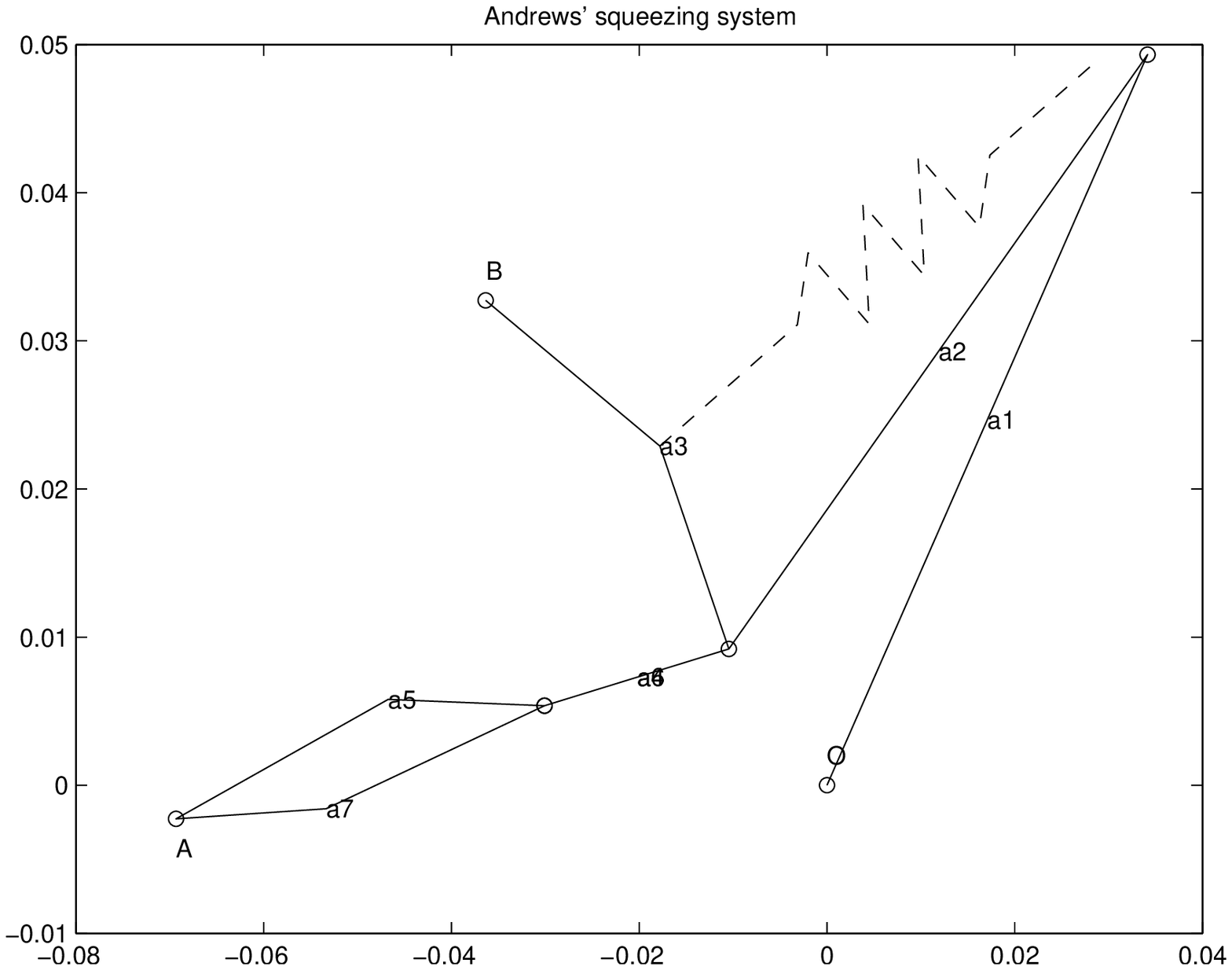, width=0.45\textwidth} \\
    \hline
  \end{tabular}
  \caption{Singular positions (according to $J_{4567}$, $T_3$) when
    $a_1=0.06,\,a_2=0.06$ and $a_3,\dots,a_7$ have the original
    values.  One can see a physical explanation to the singularity:
    the centre node $P_2$ is now 'pulled out' so that nodes $P_3$ and $P_4$ coincide.}
  \label{fig:S1_T3_sing}
\end{figure}

For example, let us concentrate on $T_5$ and choose
$a_1=0.03,\,a_2=0.055$, say, whence the system is able to reach four
singular configurations (see the left panel of Figure
\ref{fig:E_1_and_2}).  Now $s_i,c_i$ for $i=1,2,3,7$ are determined by
$\mathsf{V}(K)$.  The other values, for angles 4,5,6, are determined
by $\mathsf{V}(T_5)$.  The results are in the Table \ref{tab:J_4567}.
The corresponding configurations are visualized in Figure
\ref{fig:S1_T5_sing}.
\begin{table}[htb]
\begin{center}
\begin{tabular}{|c|c|c|c|c|}
  \hline
  variable & singularity $1$ & singularity $2$ & singularity  $3$& singularity  $4$\\ \hline
  $c_1$ & -0.8322 & -0.4564 & -0.1157 & -0.1038  \\ \hline
  $s_1$ & -0.5544 & 0.8898  & -0.9933 & 0.9946   \\ \hline
  $c_2$ & -0.3045 & -0.3045 & 0.4467  & 0.4467   \\ \hline
  $s_2$ & 0.9525  & -0.9525 & 0.8947  & -0.8947  \\ \hline
  $c_3$ & 0.4300  & 0.4300  & 0.9267  & 0.9267   \\ \hline
  $s_3$ & -0.9029 & -0.9029 & -0.3759 & -0.3759  \\ \hline
  $c_4$ & 0       & 0       & 0       & 0        \\ \hline
  $s_4$ & -1      & -1      & -1      & -1       \\ \hline
  $c_5$ & 0.0695  & 0.0695  & 0.9917  & 0.9917   \\ \hline
  $s_5$ & 0.9976  & 0.9976  & 0.1283  & 0.1283   \\ \hline
  $c_6$ & 0       & 0       & 0       & 0        \\ \hline
  $s_6$ & 1       & 1       & 1       & 1        \\ \hline
  $c_7$ & -0.9976 & -0.9976 & -0.1283 & -0.1283  \\ \hline
  $s_7$ & 0.0695  & 0.0695  & 0.9917  & 0.9917   \\ \hline
\end{tabular}
\end{center}
\vspace*{5mm}
Calculating the corresponding angles we get the following values.
\vspace*{5mm}
\begin{center}
\begin{tabular}{|c|c|c|c|c|}
  \hline
  Angle & singularity $1$ & singularity $2$ & singularity $3$ & singularity $4$\\\hline
  $y_1$ & -2.5539 & 2.0448  & -1.6867 & 1.6747   \\ \hline
  $y_2$ & 1.8802  & -1.8802 & 1.1077  & -1.1077  \\ \hline
  $y_3$ & -1.1264 & -1.1264 & -0.3853 & -0.3853  \\ \hline
  $y_4$ & -1.5708 & -1.5708 & -1.5708 & -1.5708  \\ \hline
  $y_5$ & 1.5012  & 1.5012  & 0.1287  & 0.1287   \\ \hline
  $y_6$ & 1.5708  & 1.5708  & 1.5708  & 1.5708   \\ \hline
  $y_7$ & 3.0720  & 3.0720  & 1.6995  & 1.6995   \\ \hline
\end{tabular}
\end{center}
  \caption{The singularities of $J_{4567}$ type, original values apart from $a_1,a_2$. The values are presented only with 4 decimals but were computed with 16 decimals.}
  \label{tab:J_4567}
\end{table}
Doing similar tests with $T_3$ instead of $T_5$ yields the $E_i$ areas
in the right hand panel of Figure \ref{fig:E_1_and_2}.  Singular
configurations implied by $T_3$, with choices $a_1=0.06,\, a_2=0.06$
which imply 4 singularities, are in Figure \ref{fig:S1_T3_sing}.  To
save space we have not tabulated the actual values of the angles in
$T_3$ case.

\subsection{$J_{367}$ singularity: original values, apart from $b_1,a_1,a_2$}
\label{Exam:3}

\begin{figure}[htb]
  \centering
  \epsfig{file=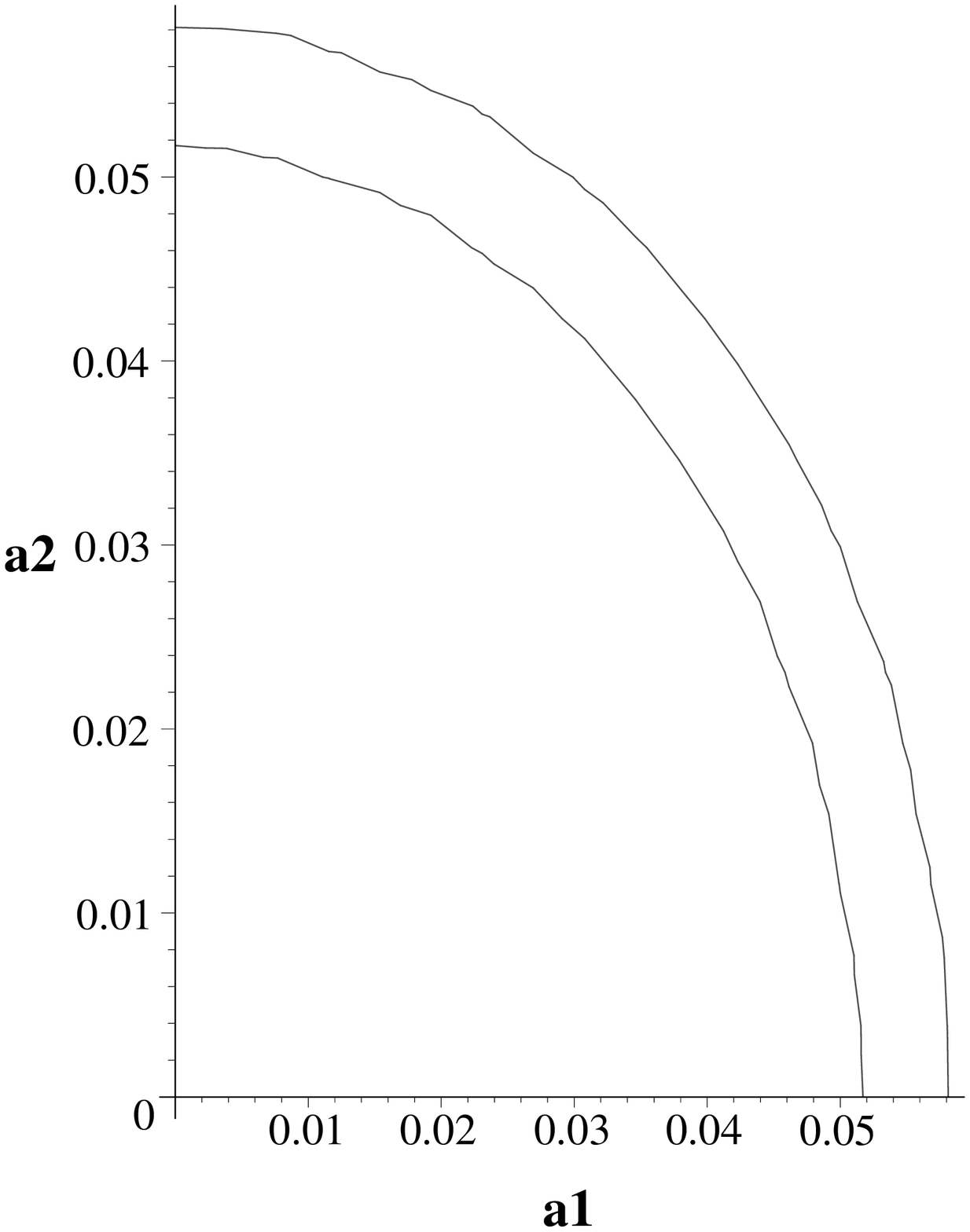, width=7cm, height=7cm}
  \caption{$J_{367},U_1$ case: the region inside the annulus is where
    $n_{3}(4a_1a_2-n_3)\ge 0$.}
  \label{fig:S2_n3n4}
\end{figure}

\begin{figure}[htb]
  \centering
  \begin{tabular}{l|r}
    \epsfig{file=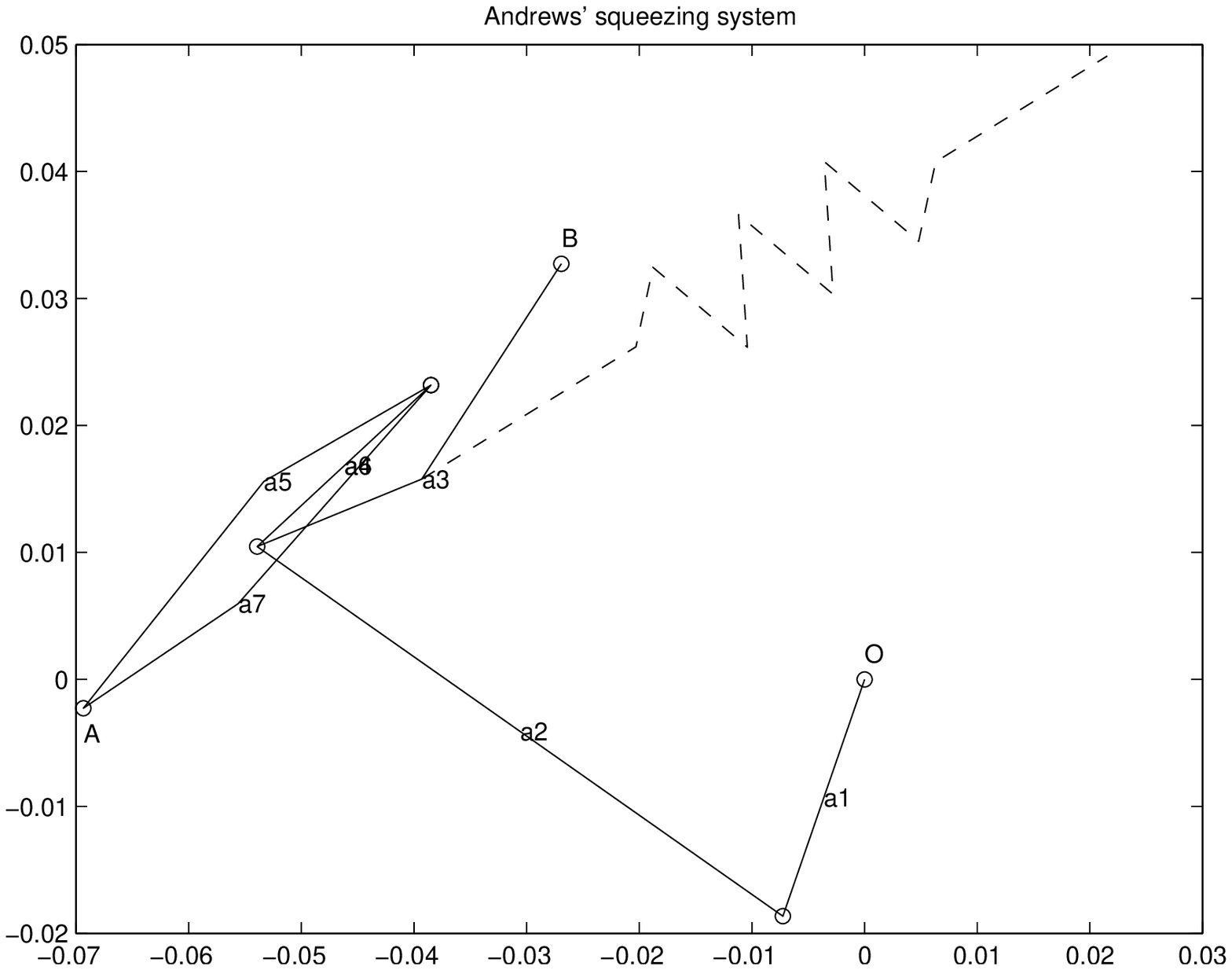, width=0.45\textwidth} &
    \epsfig{file=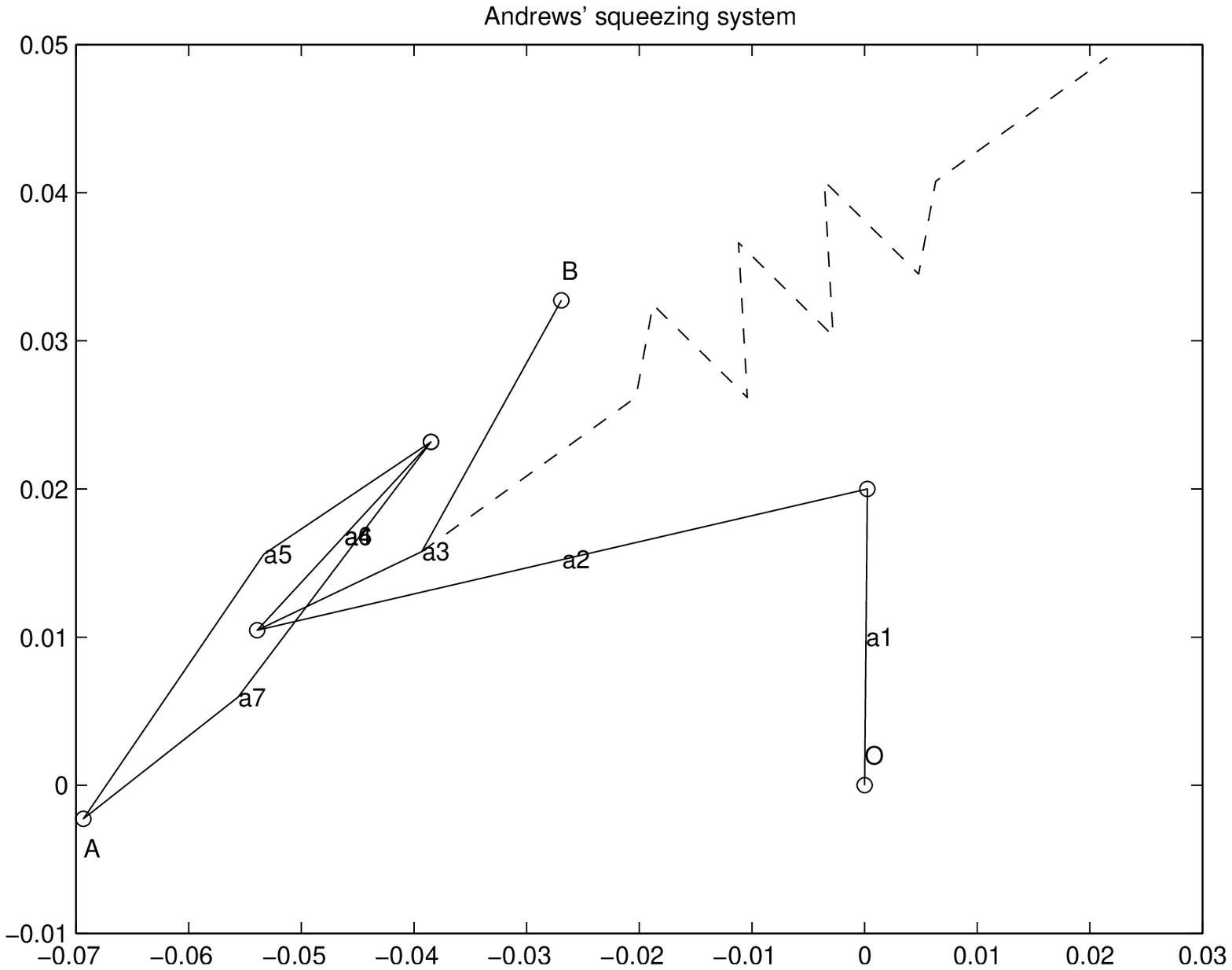, width=0.45\textwidth} \\
    \hline
  \end{tabular}
  \caption{Singular positions (according to $J_{367}$, $U_1$) when
    $b_1=-0.02691,\, a_1=0.02,\, a_2=0.055$ and $a_3,\dots,a_7$ have
    the original values.  The physical interpretation is as in
    Figure \ref{fig:S1_T5_sing}.}
  \label{fig:S2_sing}
\end{figure}

A necessary condition to have a $J_{367}$ type singularity is at least
one of the $z_i$'s vanishes \eqref{eq:U_z}.
Substituting the original parameter values we notice that none of
these is zero.  Let us then investigate how we should change some of
the parameters in order to have $J_{367}$ type singularities.  Take $b_1$
and $U_1$, say, and choose $b_1:=-0.026913593$ so that $z_1=0$. \footnote{This corresponds to  moving $B$ slightly to left.}   We seek to further
fulfil the sufficient requirements by $U_1$:
\begin{align*}
  &n_{3}(4a_1a_2-n_3)\geq 0 \tag{\ref{U1_cond_1}} \\
  & t_7 t_5 t_6 t_8\leq 0, \tag{\ref{U1_cond_2}}
\end{align*}
and use $L_1,L_2$ to find the actual singular configurations.  With
the original parameter values $t_6=0$, therefore \eqref{U1_cond_2} is
fulfilled.  Therefore we only need to study \eqref{U1_cond_1}.  For
that, we proceed analogously to Example \ref{Exam:2}: treat the
expression $n_{3}(4a_1a_2-n_3)$ as a function of $a_1,a_2$.  For that,
we first need $c_7,s_7$.  Them we get from \eqref{S2_U1}
\begin{align*}
  c_7 &= \frac{b_2-w_2}{a_6-a_3-a_7}=-0.6364 \\
  s_7 &= \frac{b_1-w_1}{a_3+a_7-a_6}=0.7714.
\end{align*}
The region of $a_1,a_2$ plane where $n_{3}(4a_1a_2-n_3)\ge 0$ is shown
in Figure \ref{fig:S2_n3n4}.
We pick a value inside the ``allowed'' annulus, say $a_1=0.02$ and
$a_2=0.055$ in order to get singularities.  Then let us find the
actual singular configurations: since $t_6=0$, from \eqref{S2_U1_c4}
we get $c_4=0$ and from \eqref{S2_U1_s4} $s_4=-1$.  The other angles
are found as follows: 3 and 6 from \eqref{S2_U1} and the remaining
ones 1,2,5 from $L$.  The results are in Table \ref{tab:J_367}.  The
corresponding singular configurations are drawn in Figure
\ref{fig:S2_sing}.  Note that there are only two singular
configurations, instead of four, since \eqref{S2_U1_c4} has only one
(double) root $c_4=0$ instead of two separate roots.
\begin{table}[htb]
\begin{center}
\begin{tabular}{|c|c|c|}
  \hline
  variables & singularity $1$ & singularity $2$\\ \hline
  $c_1$ & -0.3621 & 0.0127   \\ \hline
  $s_1$ & -0.9322 & 0.9999   \\ \hline
  $c_2$ & 0.1860  & 0.1860   \\ \hline
  $s_2$ & 0.9862  & -0.9826  \\ \hline
  $c_3$ & 0.6364  & 0.6364   \\ \hline
  $s_3$ & -0.7714 & -0.7714  \\ \hline
  $c_4$ & 0       & 0        \\ \hline
  $s_4$ & -1      & -1       \\ \hline
  $c_5$ & 0.7714  & 0.7714   \\ \hline
  $s_5$ & 0.6364  & 0.6364   \\ \hline
  $c_6$ & 0       & 0        \\ \hline
  $s_6$ & 1       & 1        \\ \hline
  $c_7$ & -0.6364 & -0.6364  \\ \hline
  $s_7$ & 0.7714  & 0.7714   \\ \hline
\end{tabular}
\end{center}
Expressed in angles, these are
\begin{center}
\begin{tabular}{|c|c|c|c|c|}
  \hline
  Angles & singularity $1$ & singularity $2$ \\\hline
  $y_1$ & -1.9413 & 1.5581   \\\hline
  $y_2$ & 1.3837  & -1.3837  \\\hline
  $y_3$ & -0.8810 & -0.8810  \\\hline
  $y_4$ & 1.5708  & 1.5708   \\\hline
  $y_5$ & 0.6898  & 0.6898   \\\hline
  $y_6$ & 1.5708  & 1.5708   \\\hline
  $y_7$ & 2.2606  & 2.2606   \\\hline
\end{tabular}
\end{center}
  \caption{The singularities of $J_{367}$ type, original values apart from $b_1,a_1,a_2$. The values are presented only with 4 decimals but were computed with 16 decimals.}
  \label{tab:J_367}
\end{table}

\subsection{A rational case}
\label{Exam:4}

\begin{figure}[htb]
  \centering
  \epsfig{file=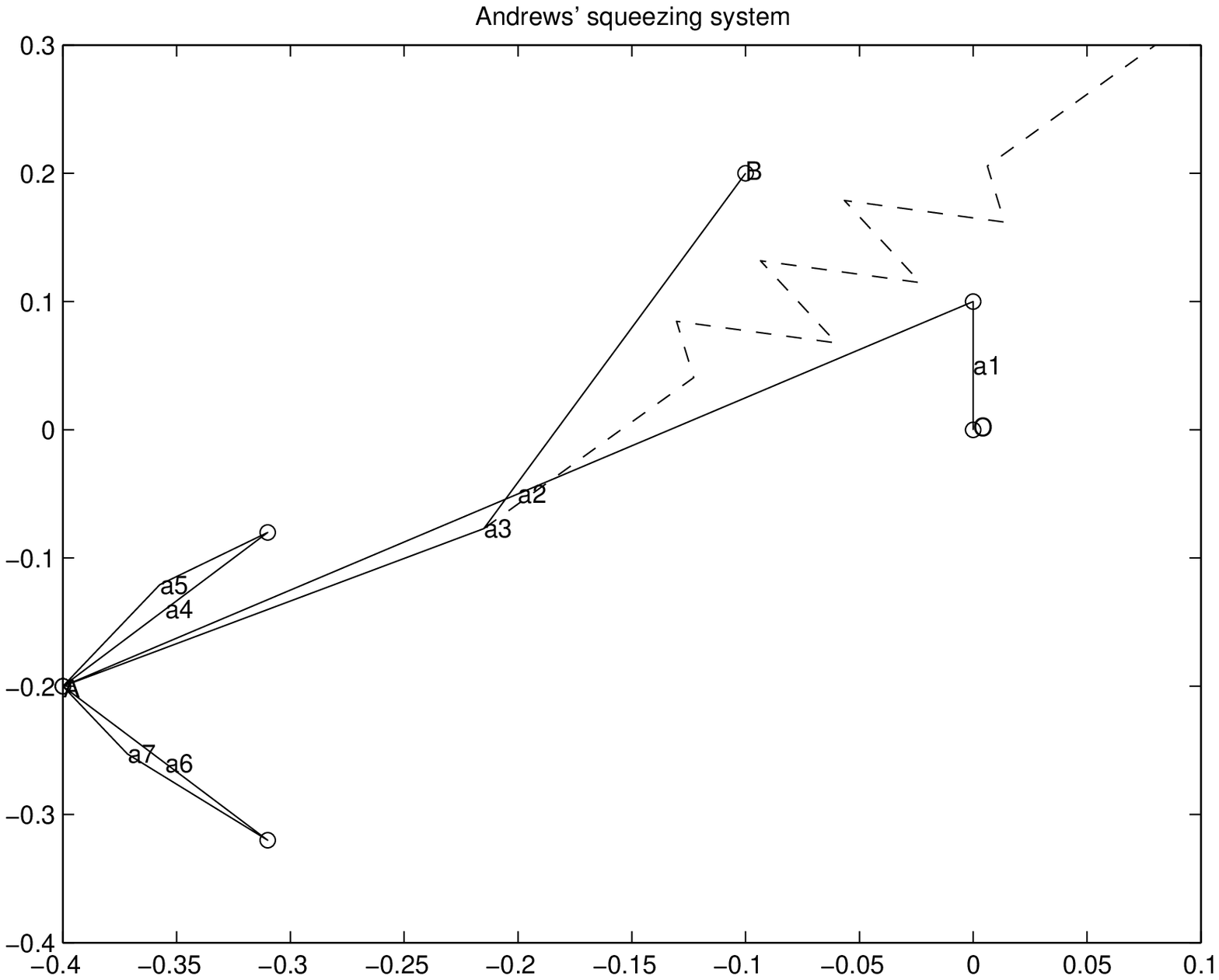, width=0.7\textwidth} 
  \caption{A singular configuration with rational $c_i,s_i,a_i,b_i$.}
  \label{fig:rat_cos_sin}
\end{figure}

Finally, let us show a rational valued singularity, that is $c_i,s_i\in\Q$.
Choose
\begin{align*}
  a_4=a_5=a_6=a_7=3/20 &\quad a_1=1/10 \quad a_2=a_3=1/2 \\
  b_1=-1/10 \quad b_2=1/5 &\quad w_1=-2/5 \quad w_2=-1/5
\end{align*}
and solve $c,s$ from the generators of $I_2 \cup J \cup F_J$ in
\eqref{eq:a4_is_a5}.  Now $c_5,s_5,c_7,s_7$ are arbitrary (apart from
$c_5^2+s_5^2=1$, $c_7^2+s_7^2=1$) and the chosen result is (see also
Figure \ref{fig:rat_cos_sin})
\begin{align*}
  c &= (0,  3/5,    4/5,  0,    3/5,     0,    4/5) \\
  s &= (1,  -4/5,   -3/5,   -1,    4/5,    1,    3/5).
\end{align*}

\section{Conclusion}
\label{sec:discuss}

We have studied singularities of the multibody system ``Andrews'
squeezing system'' which is a well-known benchmark problem both for
multibody solvers and differential-algebraic equation solvers.  
Using our tools we have shown in Theorem \ref{thm:3}
that the original benchmark problem is indeed void of singularities,
thereby assuring that whatever numerical problems in the benchmark
tests are met, they are indeed due to something else than a nearby
singularity of the system.  Apparently, this non-singularity of the
problem has not been rigorously proven in the literature.

However, we have shown that with suitably chosen parameters $(a,b,w)$, this
system can exhibit singular configurations.  In fact, there are
families of values $(a,b,w)$ that produce singularities, see Theorems
\ref{thm:1} and \ref{thm:2}. We provide examples of singularities, calculated using the original
benchmark parameter values apart from $b_1,a_1,a_2$.  Considering
$a_1,a_2$ as freely chosen parameters, Figures \ref{fig:E_1_and_2} and
\ref{fig:S2_n3n4} show the areas of $a_1,a_2$ plane where the system
exhibits singularities.  For example, choosing the point $(a_1,a_2)$
within the intersection of the three areas in Figures
\ref{fig:E_1_and_2} (both panels) and \ref{fig:S2_n3n4} would give a
system with 10 singular configurations.

A natural question that remains is, if these presented singularities are 
the only possible ones? In other words are there singularities which do not come from the singularities of some subsystem? While the Gr\"obner bases techniques {\em in principle} provide a way to answer this question directly, we could not do so in practice due to complexity problems.

\subsection{Appendix} 

\paragraph{The coefficients $f_i$:} 
The coefficients $f_1,\dots,f_5$ in the context of $T_5$ are \\
\begin{align*}
  f_1 &=4(a_5-a_4)^2(b_1^2-2b_1w_1+b_2^2-2b_2w_2+w_1^2+w_2^2)\\
  & =4(a_5-a_4)^2|b-w|^2, \\
  f_2 &= 4(w_1-b_1)(a_4-a_5)(-b_1^2+2b_1w_1-b_2^2+2b_2w_2-w_1^2-w_2^2+a_3^2-a_4^2+2a_4a_5-a_5^2)\\
  &=4(w_1-b_1)(a_4-a_5)\big(a_3^2-(a_4-a_5)^2-|b-w|^2\big), \\
  f_3 &=b_1^2-2b_1w_1+b_2^2-2b_2w_2+2b_2a_4-2b_2a_5+w_1^2+w_2^2-2w_2a_4+2w_2a_5-a_3^2+a_4^2-2a_4a_5+a_5^2\\
  &=|b-w|^2+2(b_2-w_2)(a_4-a_5) -a_3^2+(a_4-a_5)^2, \\
  f_4 &=b_1^2-2b_1w_1+b_2^2-2b_2w_2-2b_2a_4+2b_2a_5+w_1^2+w_2^2+2w_2a_4-2w_2a_5-a_3^2+a_4^2-2a_4a_5+a_5^2 \\
  &= |b-w|^2-2(b_2-w_2)(a_4-a_5) -a_3^2+(a_4-a_5)^2, \\
  f_5 &= a_3^2-a_4^2+2a_4a_5-a_5^2-b_1^2+2b_1w_1-b_2^2+2b_2w_2-w_1^2-w_2^2 \\
  &= a_3^2-(a_4-a_5)^2-|b-w|^2.
\end{align*}

\paragraph{The coefficients $d_i,l_i$:}
The coefficients $d_i,l_i$ in the context of $K_2$ are \\
\begin{tabular}[h]{ll}
  \hline\\
$ d_1 $ &= $ 2a_1a_2(a_3^2 +2a_3b_1s_3 -2a_3b_2c_3 +b_1^2+b_2^2) $ \\[3mm]
$ d_2 $ &= $ -4a_1^2a_2^2s_2^2$ \\[3mm]
$ d_3 $ &= $ -a_1^4 +2a_1^2a_2^2 +a_1^2a_3^2 +2a_1^2a_3b_1s_3 -2a_1^2a_3b_2c_3 +a_1^2b_1^2 $ \\[3mm]
  & $ +a_1^2b_2^2 -a_2^4 +a_2^2a_3^2 +2a_2^2a_3b_1s_3 -2a_2^2a_3b_2c_3 +a_2^2b_1^2 +a_2^2b_2^2 $ \\[3mm]
  \hline\\
$ l_1 $ &= $ -2a_1a_2(a_3^2 +2a_3b_1s_3 -2a_3b_2c_3 +b_1^2+b_2^2) $ \\[3mm]
$ l_2 $ &= $ 2a_1a_2(a_3s_3+b_1) $ \\[3mm]
$ l_3 $ &= $ -(a_3c_3-b_2)(a_1^2 -a_2^2 +a_3^2 +2a_3b_1s_3 -2a_3b_2c_3 +b_1^2+b_2^2) $ \\[3mm]
$ l_4 $ &= $ 2a_1a_2s_1s_2-(a_3s_3+b_1)a_2c_2+(a_3c_3-b_2)a_2s_2-(a_3s_3+b_1)a_1 $. \\[3mm]
  \hline\\
\end{tabular}

We can also simplify these expressions: 
\begin{eqnarray*}
  d_0 &=& a_3^2 +|b|^2+2a_3(b_1s_3-b_2c_3) \\
  d_1 &=& 2a_1a_2d_0 \\
  d_2 &=& n_1n_2 \\
  d_3 &=& (a_1^2+a_2^2)d_0-(a_1^2-a_2^2)^2 \\
  n_1 &=& (a_1+a_2)^2 -d_0 \\
  n_2 &=& (a_1-a_2)^2 -d_0 = 4a_1a_2-n_1 \\
  l_1 &=& -d_1 \\
  l_3 &=& -(a_3c_3-b_2)(a_1^2-a_2^2+d_0) \\
  l_4 &=& -(a_3s_3+b_1)(a_2c_2+a_1) +a_2s_2 (a_3c_3-b_2+2a_1s_1) \\
  \hat g_1 &=& -4a_1^2a_2^2s_2^2+n_1(4a_1a_2-n_1) \\
  \hat g_2 &=& 2a_1a_2 d_0 c_2 +(a_1^2+a_2^2)d_0 -(a_1^2-a_2^2)^2 \\
  \hat g_3 &=& -2a_1a_2 d_0 s_1 +2a_1a_2 (a_3 s_3+b_1) -(a_3c_3-b_2)(a_1^2-a_2^2+d_0) \\
  \hat g_4 &=& (a_1^2-a_2^2)c_1+ l_4 
\end{eqnarray*}

\paragraph{The coefficients $r_i$:} 
The coefficients $r_i$ in the context of $L_1$ are \\
\begin{tabular}[h]{ll}
  \hline\\
  $ r_1$ &= $(a_1^2+a_2^2)|b|^2-2b_1a_1^2a_3s_7-2b_1a_2^2a_3s_7+2b_2a_1^2a_3c_7+2b_2a_2^2a_3c_7-(a_1^2-a_2^2)^2 
     +(a_1^2 +a_2^2)a_3^2 $ \\[3mm]
  $r_2$ &= $2a_1(b_1a_2-a_2a_3s_7)s_2 $ \\[3mm]
  $r_3$ &= $b_1^2b_2+b_1^2a_3c_7-2b_1b_2a_3s_7-2b_1a_3^2c_7s_7+b_2^3+3b_2^2a_3c_7+b_2a_1^2-b_2a_2^2+3b_2a_3^2c_7^2 $ \\[3mm]
  & $+b_2a_3^2s_7^2+a_1^2a_3c_7-a_2^2a_3c_7+a_3^3c_7 $ \\[3mm]
  $r_4$ &= $(2a_1a_2)s_1s_2+(-b_1a_2+a_2a_3s_7)c_2+(-b_2a_2-a_2a_3c_7)s_2+(-b_1a_1+a_1a_3s_7)$\\[3mm]
  \hline\\
\end{tabular}

\end{document}